%% file: mincirc.tex
\documentclass[a4paper]{article}

\usepackage{graphicx}
\usepackage{amsmath}
\usepackage{amsthm}
\usepackage{amssymb}
\usepackage{color}
\usepackage{bbm}
\usepackage{floatflt}
\usepackage{url}

\graphicspath{{./} {./smallFigures/}}
\newcommand{\R}{\mathbbm R}
\newcommand{\C}{\mathbbm C}
\newcommand{\Z}{\mathbbm Z}
\newcommand{\cir}{\mbox{$\bigcirc\kern -7.6pt c\ $}}
\newcommand{\sph}{\mbox{$\bigcirc\kern -7.6pt s\ $}}

\theoremstyle{plain}
\newtheorem{theorem}{Theorem}
\newtheorem{proposition}{Proposition}

\newtheorem{lemma}{Lemma}

\renewcommand{\Re}{\mathop\mathrm{Re}}
\renewcommand{\Im}{\mathop\mathrm{Im}}

\theoremstyle{definition}
\newtheorem{definition}{Definition}

\theoremstyle{remark}
\newtheorem*{examples}{Examples}
\newtheorem*{remark}{Remark}
\newtheorem*{remarks}{Remarks}

\DeclareMathOperator{\im}{Im}
\DeclareMathOperator{\Li}{Li}

\author{Alexander I.~Bobenko\setcounter{footnote}{0}\thanks{Partially
    supported by the DFG Research Center ``Mathematics for Key Technologies''
    (FZT 86) in Berlin.} %
  \and Tim Hoffmann\setcounter{footnote}{2}\thanks{Supported by the DFG
    Research Center ``Mathematics for Key Technologies'' (FZT 86) in Berlin
    and the Alexander von Humboldt Foundation.}  %
  \and Boris A.~Springborn\setcounter{footnote}{3}\thanks{Supported by the
    DFG Research Center ``Mathematics for Key Technologies'' (FZT 86) in
    Berlin.}}

\title{Minimal surfaces from circle patterns:\\
  Geometry from combinatorics}

\begin{document}

\maketitle

\section{Introduction}
\label{sec:intro}

The theory of polyhedral surfaces and, more generally, the field of discrete
differential geometry are presently emerging on the border of differential
and discrete geometry. Whereas classical differential geometry investigates
smooth geometric shapes (such as surfaces), and discrete geometry studies
geometric shapes with a finite number of elements (polyhedra), the theory of
polyhedral surfaces aims at a development of discrete equivalents of the
geometric notions and methods of surface theory. The latter appears then
as a limit of the refinement of the discretization. Current progress in this
field is to a large extent stimulated by its relevance for computer graphics
and visualization.

One of the central problems of discrete differential geometry is to find
proper discrete analogues of special classes of surfaces, such as minimal,
constant mean curvature, isothermic surfaces, etc. Usually, one can suggest
various discretizations with the same continuous limit which have quite
different geometric properties. The goal of discrete differential geometry is
to find a discretization which inherits as many essential properties of the
smooth geometry as possible.

\begin{figure}[tb]%
\includegraphics[width=0.42\textwidth,clip=true]{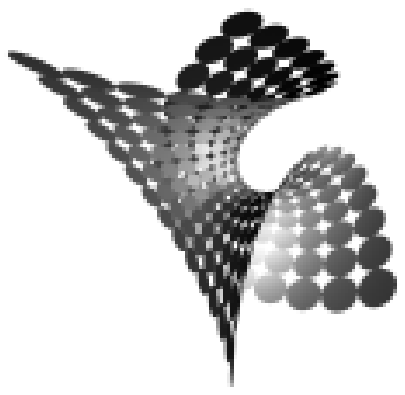}%
\includegraphics[width=0.55\textwidth,clip=true]{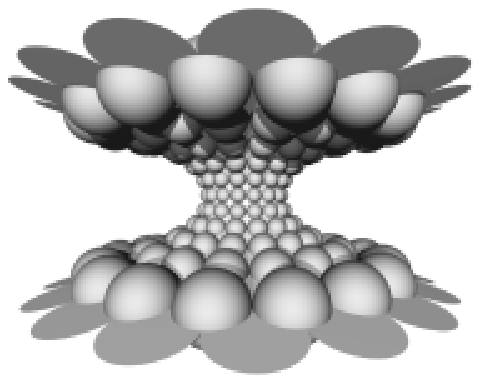}%
\caption{A discrete minimal Enneper surface {\em(left)}\/ and a discrete
  minimal catenoid {\em(right)}.}
\label{fig:enneper_catenoid}
\end{figure}

\begin{figure}[tb]%
\hfill{}%
\includegraphics[width=0.57\textwidth,clip=true]{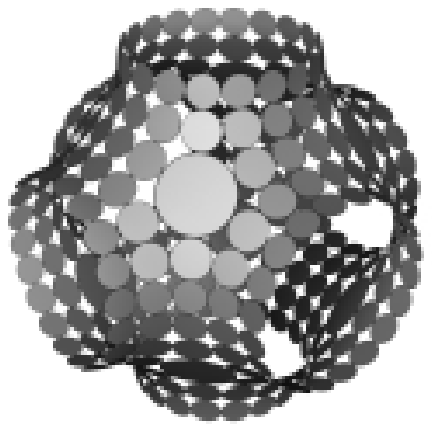}%
\hfill{}%
\includegraphics[width=0.35\textwidth,clip=true]{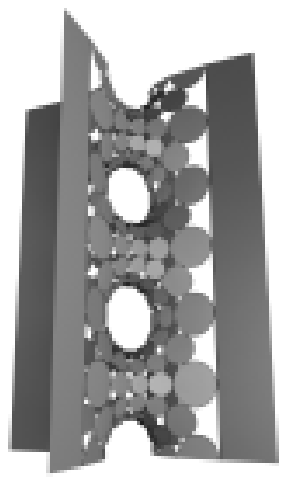}%
\hspace*{\fill}
\caption{A discrete minimal Schwarz P-surface {\em(left)}\/ and a discrete 
  minimal Scherk tower {\em(right)}.}
\label{fig:schwarz_scherk}
\end{figure}

Our discretizations are based on quadrilateral meshes, i.e.\ we discretize
parameterized surfaces. For the discretization of a special class of surfaces,
it is natural to choose an adapted parameterization. In this paper, we
investigate conformal discretizations of surfaces, i.e.\ discretizations in
terms of circles and spheres, and introduce a new discrete model for minimal
surfaces. See Figs.~\ref{fig:enneper_catenoid} and~\ref{fig:schwarz_scherk}.
In comparison with direct methods (see, in particular, \cite{PR02}), leading
usually to triangle meshes, the less intuitive discretizations of the present
paper have essential advantages: they respect conformal properties of
surfaces, possess a maximum principle (see Remark on
p.~\pageref{rem:maximum_principle}), etc.

We consider minimal surfaces as a subclass of isothermic surfaces. The
analogous discrete surfaces, {\em discrete S-isothermic surfaces} \cite{BP99}
consist of touching spheres, and of circles which intersect the spheres
orthogonally in their points of contact. See Fig.~\ref{fig:enneper_catenoid}
{\em(right)}. Continuous isothermic surfaces allow a duality transformation,
the Christoffel transformation. Minimal surfaces are characterized among
isothermic surfaces by the property that they are dual to their Gauss map.
The duality transformation and the characterization of minimal surfaces
carries over to the discrete domain. Thus, one arrives at the notion of {\em
  discrete minimal S-isothermic surfaces}, or {\em discrete minimal surfaces}
for short. The role of the Gauss maps is played by discrete S-isothermic
surfaces, the spheres of which all intersect one fixed sphere orthogonally.
Due to a classical theorem of Koebe (see Section~\ref{sec:koebe}) any
3-dimensional combinatorial convex polytope can be (essentially uniquely)
realized as such a Gauss map.

This definition of discrete minimal surfaces leads to a construction method
for discrete S-isothermic minimal surfaces from discrete holomorphic data, a
form of a discrete Weierstrass representation (see
Section~\ref{sec:weierstrass}). Moreover, the classical ``associated family''
of a minimal surface, which is a one-parameter family of isometric
deformations preserving the Gauss map, carries over to the discrete setup
(see Section~\ref{sec:assoc_family}).

Our general method to construct discrete minimal surfaces is
schematically shown in the following diagram. (See also
Fig.\ref{fig:recipe}.)
\begin{center}
  continuous minimal surface \\
  {\large$\Downarrow$} \\
  image of curvature lines under Gauss-map\\
  {\large$\Downarrow$} \\
  cell decomposition of (a branched cover of) the sphere\\
  {\large$\Downarrow$} \\
  orthogonal circle pattern \\
  {\large$\Downarrow$} \\
  Koebe polyhedron\\
  {\large$\Downarrow$} \\
  discrete minimal surface
\end{center}
As usual in the theory on minimal surfaces \cite{HK97}, one starts
constructing such a surface with a rough idea of how it should look. To use
our method, one should understand its Gauss map and the {\em combinatorics}\/
of the curvature line pattern. The image of the curvature line pattern under
the Gauss map provides us with a cell decomposition of (a part of) $S^2$ or a
covering.  From these data, applying the Koebe theorem, we obtain a circle
packing with the prescribed combinatorics.  Finally, a simple dualization
step yields the desired discrete minimal surface.

Let us emphasize that our data, besides possibly boundary conditions, are
purely combinatorial---the combinatorics of the curvature line pattern. All
faces are quadrilaterals and typical vertices have four edges. There may exist
distinguished vertices (corresponding to the ends or umbilic points of a
minimal surface) with a different number of edges.

The most nontrivial step in the above construction is the third one listed in
the diagram. It is based on the Koebe theorem. It implies the existence and
uniqueness for the discrete minimal S-isothermic surface under consideration,
but not only this. This theorem can be made an effective tool in constructing
these surfaces. For that purpose, we use a variational principle from
\cite{BS02}, \cite{SprPhD} for constructing circle patterns. This principle
provides us with a variational description of discrete minimal S-isothermic
surfaces and makes possible a solution of some Plateau problems as well.

In Section~\ref{sec:convergence}, we prove the convergence of
discrete minimal S-isothermic surfaces to smooth minimal surfaces. The proof
is based on Schramm's approximation result for circle patterns with the
combinatorics of the square grid~\cite{Sch97}. The best known convergence
result for circle patterns is $C^{\infty}$-convergence of circle
packings~\cite{HS98}. It is an interesting question whether the convergence
of discrete minimal surfaces is as good.

Because of the convergence, the theory developed in this paper may be used to
obtain new results in the theory of smooth minimal surfaces. A typical
problem in the theory of minimal surfaces is to decide whether surfaces with
some required geometric properties exist, and to construct them. The
discovery of the Costa-Hoffman-Meeks surface \cite{HM85}, a turning point of
the modern theory of minimal surfaces, was based on the Weierstrass
representation. This powerful method allows the construction of important
examples. On the other hand, it requires a specific study for each example;
and it is difficult to control the embeddedness. Kapouleas \cite{Kap97}
proved the existence of new embedded examples using a new method. He
considers finitely many catenoids with the same axis and planes orthogonal to
this axis and shows that one can desingularize the circles of intersection by
deformed Scherk towers. This existence result is very intuitive, but it gives
no lower bound for the genus of the surfaces.  Although some examples with
lower genus are known (the Costa-Hoffman-Meeks surface and
generalizations~\cite{HM90}), proving the existence of Kapouleas' surfaces
with given genus and to construct them using conventional methods is very
difficult~\cite{WW02}. Our method may become helpful to address these
problems. At the present time, however, the construction of new minimal
surfaces from discrete ones remains a challenge.

Apart from discrete minimal surfaces, there are other interesting subclasses
of S-isothermic surfaces. In future publications, we plan to treat discrete
constant mean curvature surfaces in Euclidean 3-space and Bryant surfaces
\cite{Bry87}, \cite{CHR01}. (Bryant surfaces are surfaces with constant mean
curvature $1$ in hyperbolic 3-space.) Both are special subclasses of
isothermic surfaces that can be characterized in terms of surface
transformations.  (See~\cite{BP99} and~\cite{HHP99} for a definition of
discrete constant mean curvature surfaces in $\R^3$ in terms of
transformations of isothermic surfaces. See~\cite{HMN01} for the
characterization of continuous Bryant surfaces in terms of surface
transformations.)

More generally, we believe that the classes of
discrete surfaces considered in this paper will be helpful in the development
of a theory of discrete conformally parameterized surfaces.

\section{Discrete S-isothermic surfaces}
\label{sec:discrete_S-isothermic}


Every smooth immersed surface in 3-space admits curvature line parameters
away from umbilic points, and every smooth immersed surface admits conformal
parameters. But not every surface admits a curvature line parameterization
that is at the same time conformal.

\begin{definition}
\label{def:isothermic} A smooth immersed surface in $\R^3$ is called {\em
  isothermic}\/ if it admits a conformal curvature line parameterization in
  a neighborhood of every non-umbilic point.
\end{definition}


Geometrically, this means that the curvature lines divide an isothermic
surface into infinitesimal squares. An isothermic immersion (a surface patch
in conformal curvature line parameters)
\begin{eqnarray*}
f:{\R}^2\supset D&\to& {\R}^3\\
 (x,y)  &\mapsto& f(x,y)
\end{eqnarray*}
is characterized by the properties
\begin{equation}
\label{eq:isothermic} \|f_x\|=\|f_y\|, \ f_x\bot f_y,\ f_{xy}\in \
{\rm span}\{f_x, f_y\}.
\end{equation}
Being an isothermic surface is a M{\"o}bius-invariant property: A M{\"o}bius
transformation of Euclidean 3-space maps isothermic surfaces to isothermic
surfaces. The class of isothermic surfaces contains all surfaces of
revolution, all quadrics, all constant mean curvature surfaces, and, in
particular, all minimal surfaces (see Theorem~\ref{thm:christoffel}). In this
paper, we are going to find a discrete version of minimal surfaces by
characterizing them as a special subclass of isothermic surfaces (see Section
\ref{sec:discrete_minimal}). 

While the set of umbilic points of an isothermic surface can in general be
more complicated, we are only interested in surfaces with isolated umbilic
points, and also in surfaces all points of which are umbilic. In the case of
isolated umbilic points, there are exactly two orthogonally intersecting
curvature lines through every non-umbilic point. An umbilic point has an even
number $2k$ $(k\neq 2)$\/ of curvature lines originating from it, evenly
spaced at $\pi/k$\/ angles. Minimal surfaces have isolated umbilic points.
If, on the other hand, every point of the surface is umbilic, then the
surface is part of a sphere (or plane) and every conformal parameterization is
also a curvature line parameterization.

Definition~\ref{def:discrete_isothermic} of discrete isothermic surfaces was
already suggested in \cite{BP96}. Roughly speaking, a discrete isothermic
surface is a polyhedral surface in 3-space all faces of which are conformal
squares. To make this more precise, we use the notion of a ``quad-graph'' to
describe the combinatorics of a discrete isothermic surface, and we define
``conformal square'' in terms of the cross-ratio of four points in $\R^3$.

A cell decomposition $\mathcal D$ of an oriented two-dimensional manifold
(possibly with boundary) is called a {\em quad-graph}, if all its faces are
quadrilaterals, that is, if they have four edges. The cross-ratio of four
points $z_1$, $z_2$, $z_3$, $z_4$ in the Riemann sphere
$\widehat\C=\C\cup\{\infty\}$ is 
\begin{equation*}
  cr(z_1, z_2, z_3, z_4)=\frac{(z_1-z_2)(z_3-z_4)}{(z_2-z_3)(z_4-z_1)}.
\end{equation*}
The {\em cross-ratio of four points in\/ $\R^3$} can be defined as follows:
Let $S$ be a sphere (or plane) containing the four points. $S$ is unique
except when the four points lie on a circle (or line). Choose an orientation
on $S$ and an orientation preserving conformal map from $S$ to the Riemann
sphere.  The cross-ratio of the four points in $\R^3$ is defined as the
cross-ratio of the four images in the Riemann sphere. The two orientations on
$S$ lead to complex conjugate cross-ratios. Otherwise, the cross-ratio does
not depend on the choices involved in the definition: neither on the
conformal map to the Riemann sphere, nor on the choice of $S$ when the four
points lie in a circle. The cross-ratio of four points in $\R^3$ is thus
defined up to complex conjugation. (For an equivalent definition involving
quaternions, see \cite{BP96}, \cite{HJ03}.) The cross-ratio of four points in
$\R^3$ is invariant under M{\"o}bius transformations of $\R^3$. Conversely, if
$p_1$, $p_2$, $p_3$, $p_4\in\R^3$ have the same cross-ratio (up to complex
conjugation) as $p_1'$, $p_2'$, $p_3'$, $p_4'\in\R^3$ , then there is a
M{\"o}bius transformation of $\R^3$ which maps each $p_j$ to $p_j'$.

\begin{figure}
\hfill%
\input{conformalsquare.pstex_t}%
\hfill%
\includegraphics{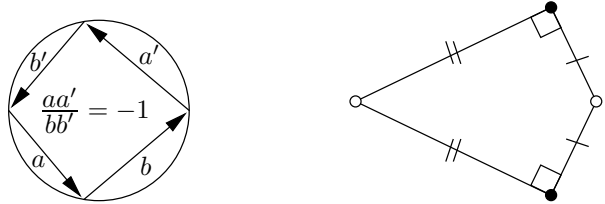}%
\hspace*{\fill}%
\caption{{\em Left:} A conformal square. The sides $a$, $a'$, $b$,
  $b'$ are interpreted as complex numbers. {\em Right:} Right-angled kites
  are conformal squares.}
\label{fig:conformalsquare}
\end{figure}

Four points in $\R^3$ form a {\em conformal square}, if their cross-ratio is
$-1$, that is, if they are M{\"o}bius-equivalent to a square.  The points of a
conformal square lie on a circle (see Fig.~\ref{fig:conformalsquare}).

\begin{definition}\label{def:discrete_isothermic}
  Let $\mathcal D$ be a quad-graph such that the degree of every interior
  vertex is even.  (That is, every vertex has an even number of edges.) Let
  $V(\mathcal D)$ be the set of vertices of $\mathcal D$. A function
  \begin{equation*}
    f:V({\mathcal D})\to\R^3
  \end{equation*}
  is called a {\em discrete isothermic surface}, if for every face of $\mathcal
  D$ with vertices $v_1$, $v_2$, $v_3$, $v_4$ in cyclic order, the points
  $f(v_1)$, $f(v_2)$, $f(v_3)$, $f(v_4)$ form a conformal square.
\end{definition}

The following three points should motivate this definition. 

\begin{itemize}
\item Like the definition of isothermic surfaces, this definition of discrete
  isothermic surfaces is M\"obius-invariant.
\item If $f:\R^2\supset D \rightarrow\R^3$ is an immersion, then for
  $\epsilon\rightarrow 0$,
  \begin{equation*}
    cr\big(f(x-\epsilon,y-\epsilon), f(x+\epsilon,y-\epsilon), 
    f(x+\epsilon,y+\epsilon),f(x-\epsilon,y+\epsilon)\big)
    =-1+O(\epsilon^2)
  \end{equation*}
  for all $x\in D$\/ if and only if $f$ is an isothermic immersion (see
  \cite{BP96}).
\item The Christoffel transformation, which also characterizes isothermic
  surfaces, has a natural discrete analogue (see Propositions~\ref{prop:dual}
  and~\ref{prop:discrete_dual}). The condition that all vertex degrees have
  to be even is used in Proposition~\ref{prop:discrete_dual}.
\end{itemize}




Interior vertices with degree different from $4$ play the role of umbilic
points. At all other vertices, two edge paths---playing the role of curvature
lines---intersect transversally. It is tempting to visualize a discrete
isothermic surface as a polyhedral surface with planar quadrilateral faces.
However, one should keep in mind that those planar faces are not M\"obius
invariant. On the other hand, when we will define discrete minimal surfaces
as special discrete isothermic surfaces, it will be completely legitimate to
view them as polyhedral surfaces with planar faces because the class of
discrete minimal surfaces is not M\"obius invariant anyway.

The Christoffel transformation~\cite{Chr1867} (see~\cite{HJ03} for a modern
treatment) transforms an isothermic surface into a dual isothermic surface.
It plays a crucial role in our considerations. For the reader's convenience,
we provide a short proof of Proposition~\ref{prop:dual}.

\begin{proposition}\label{prop:dual} 
  Let $f:{\R}^2\supset D\to {\mathbb R}^3$ be an isothermic immersion, where
  $D$\/ is simply connected. Then the formulas
  \begin{equation}
    \label{eq:dual}
    f^*_x=\dfrac{f_x}{\|f_x\|^2},\quad f^*_y=-\dfrac{f_y}{\|f_y\|^2}.  
  \end{equation}
  define (up to a translation) another isothermic immersion $f^*:{\mathbb
    R}^2\supset D\to {\mathbb R}^3$ which is called the {\em Christoffel
    transformed} or {\em dual isothermic surface}.
\end{proposition}
\begin{proof}
  First, we need to show that the 1-form $df^*=f^*_x\,dx + f^*_y\,dy$\/ is
  closed and thus defines an immersion $f^*$. From
  equations~\eqref{eq:isothermic}, we have $f_{xy}=a f_x + b f_y$, where $a$
  and $b$ are functions of $x$ and $y$. Taking the derivative of
  equations~\eqref{eq:dual} with respect to $y$ and $x$, respectively, we
  obtain
  \begin{equation*}
    f^*_{xy}=\frac{1}{\|f_x\|^2}(-a f_x + b f_y) 
    = -\frac{1}{\|f_y\|^2}(a f_x - b f_y) = f^*_{yx}.
  \end{equation*}
  Hence, $df^*$ is closed. Obviously, $\|f^*_x\|=\|f^*_y\|$, $f^*_x\bot
  f^*_y$, and $f^*_{xy}\in \operatorname{span}\{f^*_x, f^*_y\}$. Hence, $f^*$
  is isothermic.
\end{proof}

\begin{remarks}
  {\em (i)}\/ In fact, the Christoffel transformation characterizes isothermic
  surfaces: If $f$ is an immersion and equations~\eqref{eq:dual} do define
  another surface, then $f$ is isothermic.
  
  {\em (ii)}\/ The Christoffel transformation is not M\"obius invariant: The
  dual of a M\"obius transformed isothermic surface is not a M\"obius
  transformed dual.
  
  {\em (iii)}\/ In equations~\eqref{eq:dual}, there is a minus sign in the
  equation for $f^*_y$, but not in the equation for $f^*_x$. This is an
  arbitrary choice. Also, a different choice of conformal curvature line
  parameters, this means choosing $(\lambda x,\lambda y)$ instead of $(x,y)$,
  leads to a scaled dual immersion. Therefore, it makes sense to consider the
  dual isothermic surface as defined only up to translation and (positive or
  negative) scale.
\end{remarks}

The Christoffel transformation has a natural analogue in the discrete
setting: In Proposition~\ref{prop:discrete_dual}, we define the dual discrete
isothermic surface. The basis for the discrete construction is the following
Lemma. Its proof is straightforward algebra.

\begin{lemma}
  \label{le:discrete_dual}
  Suppose $a, b, a', b'\in\mathbb{C}\setminus\{0\}$ with
  \[a + b + a' + b' = 0, \quad \frac{aa'}{bb'}=-1\]
  and let
  \[a^* = \frac{1}{\,\overline{a}\,},\quad
  {a'}^* = \frac{1}{\,\overline{a}'\,},\quad b^* =
  -\frac{1}{\,\overline{b}\,},\quad {b'}^* = -\frac{1}{\,\overline{b'}\,},\]
  where $\overline{z}$ denotes the complex conjugate of $z$.  Then
  \[a^* + b^* + {a'}^* + {b'}^* = 0, \quad \frac{a^*{a'}^*}{b^*{b'}^*}=-1.\]
\end{lemma}

\begin{figure}%
\centering%
\includegraphics[width=0.3\textwidth]{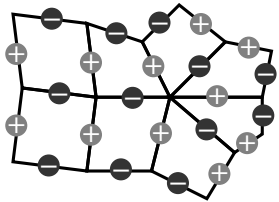}%
\caption{Edge labels of a discrete isothermic surface.}%
\label{fig:plusminusedges}%
\end{figure}

\begin{proposition}\label{prop:discrete_dual}
  Let $f:V({\mathcal D})\rightarrow\R^{3}$ be a discrete isothermic surface,
  where the quad-graph $\mathcal D$ is simply connected. Then the edges of $\mathcal
  D$ may be labelled ``$\,+$'' and ``$\,-$'' such that each quadrilateral has
  two opposite edges labelled ``$\,+$'' and the other two opposite edges
  labeled ``$\,-$'' (see Fig.~\ref{fig:plusminusedges}). The {\em dual
    discrete isothermic surface}\/ is defined by the formula
  \begin{equation*}
    \Delta f^*=\pm \dfrac{\Delta f}{\|\Delta f\|^2},
  \end{equation*}
  where $\Delta f$ denotes the difference of neighboring vertices and the
  sign is chosen according to the edge label.
\end{proposition}

For a consistent edge labelling to be possible, it is necessary that each
vertex has an even number of edges. This condition is also sufficient if the
the surface is simply connected.


In Definition~\ref{def:S-quad-graph} we define S-quad-graphs. These are
specially labeled quad-graphs that are used in
Definition~\ref{def:S-isothermic} of S-isothermic surfaces. These
S-isothermic surfaces form the subclass of discrete isothermic surfaces that we
use to define discrete minimal surfaces in
Section~\ref{sec:discrete_minimal}.  For a discussion why S-isothermic
surfaces are the right class to consider, see the Remark at the end of
Section~\ref{sec:discrete_minimal}.


\begin{definition}
\label{def:S-quad-graph} 
An {\em S-quad-graph} is a quad-graph $\mathcal D$\/ with interior vertices
of even degree as in Definition~\ref{def:discrete_isothermic} and the
following additional properties (see Fig.~\ref{fig:OneTwo+S_quad-graph}):
\begin{enumerate}
\renewcommand{\theenumi}{{\em(\roman{enumi})}}
\renewcommand{\labelenumi}{\theenumi}
\item The 1-skeleton of $\mathcal D$\/ is bipartite and the vertices
  are bi-colored ``black'' and ``white''. (Then each quadrilateral has two
  black vertices and two white vertices.)
\item Interior black vertices have degree $4$.
\item The white vertices be labeled $\cir$ and $\sph$ in such a
  way that each quadrilateral has one white vertex labeled $\cir$ and one
  white vertex labeled $\sph$.
\end{enumerate}
\end{definition}

\begin{definition}
\label{def:S-isothermic}
Let $\mathcal D$ be an S-quad-graph, and let $V_{b}({\mathcal D})$ be the set
of black vertices. A {\em discrete S-isothermic surface} is a map
\begin{equation*}
  f_{b}:V_{b}({\mathcal D})\rightarrow\mathbb{R}^3,
\end{equation*}
with the following properties:
\begin{enumerate}
\renewcommand{\theenumi}{{\em(\roman{enumi})}}
\renewcommand{\labelenumi}{\theenumi}
\item If $v_1, \ldots, v_{2n}\in V_{b}({\mathcal D})$ are the neighbors of a
  $\cir$-labeled vertex in cyclic order, then $f_{b}(v_1), \ldots,
  f_{b}(v_{2n})$ lie on a circle in $\mathbb{R}^3$ in the same cyclic order.
  This defines a map from the $\cir$-labeled vertices to the set of circles
  in $\mathbb{R}^3$.
\item If $v_1, \ldots, v_{2n}\in V_{b}({\mathcal D})$ are the neighbors of an
  $\sph$-labeled vertex, then $f_{b}(v_1), \ldots, f_{b}(v_{2n})$ lie on a
  sphere in $\mathbb{R}^3$. This defines a map from the $\sph$-labeled
  vertices to the set of spheres in $\mathbb{R}^3$.
\item If $v_c$ and $v_s$ are the $\cir$-labeled and the $\sph$-labeled vertex
  of a quadrilateral of $\mathcal D$, then the circle corresponding to $v_c$
  intersects the sphere corresponding to $v_s$ orthogonally. 
\end{enumerate}
\end{definition}

There are two spheres through each black vertex, and the orthogonality
condition~{\em(iii)}\/ of Definition~\ref{def:S-isothermic} implies that they
touch.  Likewise, the two circles at a black vertex touch, i.e.~they have a
common tangent at the single point of intersection. Discrete S-isothermic
surfaces are therefore composed of touching spheres and touching circles with
spheres and circles intersecting orthogonally. Interior white vertices of
degree unequal to $4$ are analogous to umbilic points of smooth isothermic
surfaces.  Generically, the orthogonality condition {\em (iii)}\/ follows
from the seemingly weaker condition that the two circles through a black
vertex touch:

\begin{figure}%
\hfill%
\includegraphics[width=0.32\textwidth]{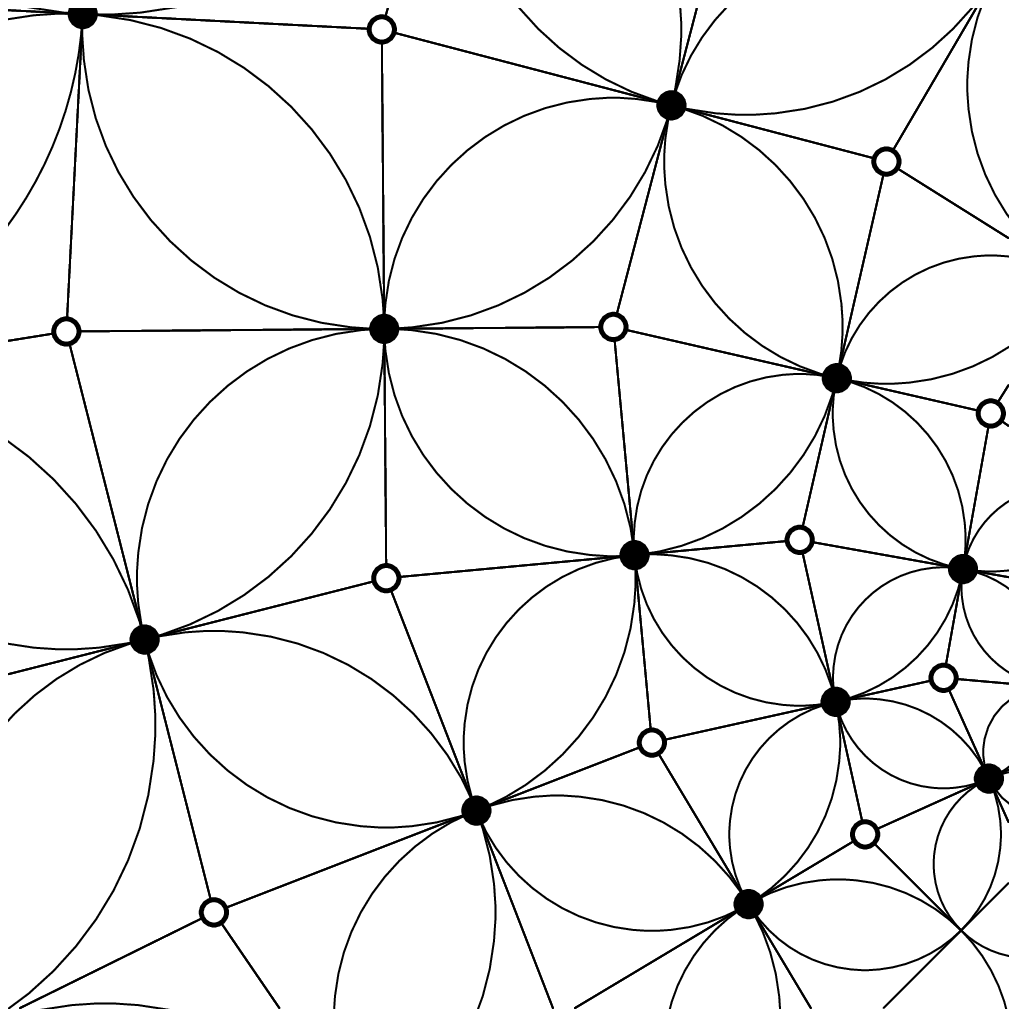}%
\hfill%
\includegraphics[width=0.4\textwidth]{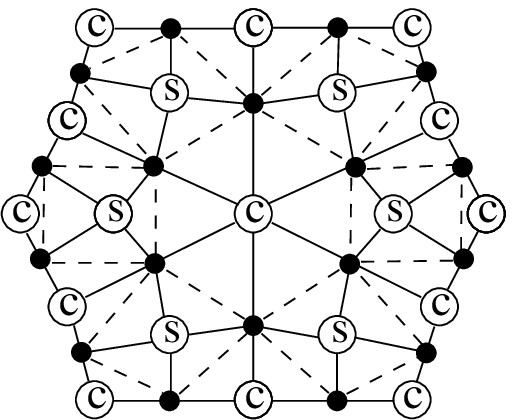}%
\hspace*{\fill}%
\caption{{\em Left:}~Schramm's circle patterns as discrete conformal
  maps. {\em Right:}~The combinatorics of S-quad-graphs.}%
\label{fig:OneTwo+S_quad-graph}%
\end{figure}

\begin{lemma} {\bf (Touching Coins Lemma)}
\label{lem:touching_coins} Whenever four circles in 3-space
touch cyclically but do not lie on a common sphere, they intersect the sphere
which passes through the points of contact orthogonally.
\end{lemma}

From any discrete S-isothermic surface, one obtains a discrete isothermic
surface (as in Definition~\ref{def:discrete_isothermic}) by adding the
centers of the spheres and circles:

\begin{definition}\label{def:central_extension}
  Let $f_{b}:V_{b}({\mathcal D})\rightarrow\mathbb{R}^3$ be a discrete
  S-isothermic surface. The {\em central extension} of $f_b$ is the discrete
  isothermic surface $f:V\rightarrow\mathbb{R}^3$ defined by
  \begin{equation*}
    f(v) = f_b(v)\quad\text{if $v\in V_b$},
  \end{equation*}
  and otherwise
  \begin{equation*}
    f(v) = \text{the center of the circle or sphere corresponding to $v$}.
  \end{equation*}  
\end{definition}

The central extension of a discrete S-isothermic surface is indeed a discrete
isothermic surface: The quadrilaterals corresponding to the faces of the
quad-graph are planar right-angled kites (see Fig.~\ref{fig:conformalsquare}
{\em(right)}) and therefore conformal squares. 

The following statement is easy to see \cite{BP99}. It says that the duality
transformation preserves the class of discrete S-isothermic surfaces.
\begin{proposition}
  \label{pro:S-dual} The Christoffel dual of a central extension of a
  discrete S-isothermic surface is a central extension of a discrete
  S-isothermic surface.
\end{proposition}

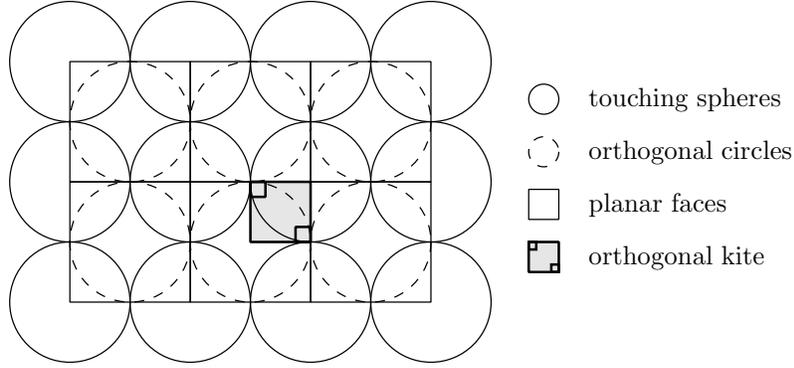
\begin{figure}[tbp]
\centerline{\input{sisothermic.pstex_t}} \caption{Geometry of a
discrete S-isothermic surface without ``umbilics''.}
\label{fig:S-isothermic}
\end{figure}

The construction of the central extension does depend on the choice of a
point at infinity, because the centers of circles and spheres are not
M\"obius invariant. Strictly speaking, a discrete S-isothermic surface has a
$3$-parameter family of central extensions. However, we will assume that one
infinite point is chosen once and for all and we will not distinguish between
S-isothermic surfaces and their central extension. Then it also makes sense
to consider the S-isothermic surfaces as polyhedral surfaces. Note that all
planar kites around a $\cir$-labeled vertex lie in the same plane: the plane
that contains the corresponding circle. We will therefore consider an
S-isothermic surface as a polyhedral surface whose faces correspond to
$\cir$-labeled vertices of the quad-graph, whose vertices correspond to
$\sph$-labeled vertices of the quad-graph, and whose edges correspond to the
black vertices of the quad-graph. The elements of a discrete S-isothermic
surface are shown schematically in Fig.~\ref{fig:S-isothermic}.
Hence:

{\em A discrete S-isothermic surface is a polyhedral surface such that the
  faces have inscribed circles and the inscribed circles of neighboring faces
  touch their common edge in the same point.}

In view of the Touching Coins Lemma (Lemma~\ref{lem:touching_coins}), this
could almost be an alternative definition.

The following lemma, which follows directly from
Lemma~\ref{le:discrete_dual}, describes the dual discrete S-isothermic
surface in terms of the corresponding polyhedral discrete S-isothermic
surface.
\begin{lemma}
\label{le:S-dual} Let $P$ be a planar polygon with an even number
of cyclically ordered edges given by the vectors
$l_1,\ldots,l_{2n}\in{\mathbb R}^2$, $l_1+\ldots+l_{2n}=0$. Suppose the
polygon has an inscribed circle with radius $R$. Let $r_j$ be the distances
from the vertices of $P$ to the nearest touching point on the circle:\/ $\|
l_j\|=r_j+r_{j+1}$. Then the vectors $l_1^*,\ldots,l_{2n}^*$ given by
$$
l_{j}^*=(-1)^{j}\dfrac{1}{r_j r_{j+1}} l_j
$$
form a planar polygon with an inscribed circle with radius\/ $1/R$.
\end{lemma}

It follows that the radii of corresponding spheres and circles of a discrete
S-isothermic surface and its dual are reciprocal.

\section{Koebe polyhedra}   \label{sec:koebe}

In this section we construct special discrete S-isothermic
surfaces, which we call the Koebe polyhedra, coming from circle
packings (and more general orthogonal circle patterns) in $S^2$.

\begin{figure}[tbp]
\hfill%
\includegraphics[width=0.3\textwidth]{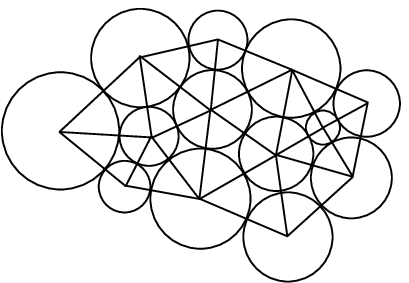}%
\hfill%
\includegraphics[width=0.3\textwidth]{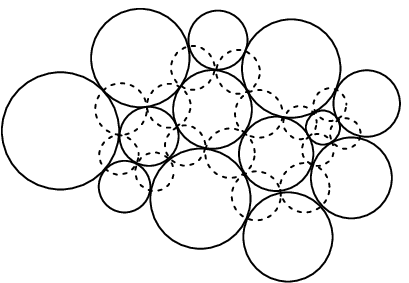}%
\hfill%
\includegraphics[width=0.3\textwidth]{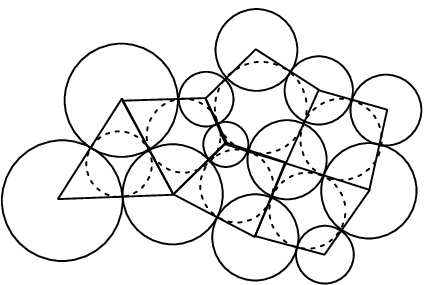}%
\hspace*{\fill}%
\caption{{\em Left:} A circle packing corresponding to a triangulation.  {\em
    Middle:} The orthogonal circles. {\em Right:} A circle packing
    corresponding to a cellular decomposition with orthogonal circles.}
 \label{fig:Koebe-combinatorics}
\end{figure}

A {\em circle packing} in $S^2$ is a configuration of disjoint
discs which may touch but not intersect. Associating vertices to
the discs and connecting the vertices of touching discs by edges
one obtains a combinatorial representation of a circle packing,
see Fig.~\ref{fig:Koebe-combinatorics} {\em(left)}.

In 1936, Koebe published the following remarkable statement about circle
packings in the sphere \cite{Koe36}.

\begin{theorem}{\bf (Koebe)}
\label{th:koebe}
 For every triangulation of the sphere there is a
packing of circles in the sphere such that circles correspond to
vertices, and two circles touch if and only if the corresponding
vertices are adjacent. This circle packing is unique up to
M{\"o}bius transformations of the sphere.
\end{theorem}

\begin{floatingfigure}{0.42\textwidth}%
\centering%
\includegraphics[width=0.2\textwidth]{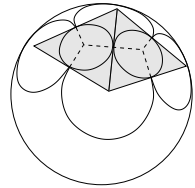}%
\caption{The Koebe polyhedron as a discrete S-isothermic surface.}%
\label{fig:Koebe-polyhedron}%
\end{floatingfigure}

Observe that for a triangulation one automatically obtains not one but two
orthogonally intersecting circle packings as shown in
Fig.~\ref{fig:Koebe-combinatorics} {\em(middle)}. Indeed, the circles passing
through the points of contact of three mutually touching circles intersect
these orthogonally. This observation leads to the following generalization of
Koebe's theorem to cellular decompositions of the sphere with faces which are
not necessarily triangular, see Fig.~\ref{fig:Koebe-combinatorics}
{\em(right)}.

\begin{theorem}
\label{thm:orthoKoebe} For every polytopal\footnote{We call a cellular
  decomposition of a surface {\em polytopal}, if the closed cells are closed
  discs, and two closed cells intersect in one closed cell if at all.}
cellular decomposition of the sphere, there exists a pattern of circles in
the sphere with the following properties.  There is a circle corresponding to
each face and to each vertex.  The vertex circles form a packing with two
circles touching if and only if the corresponding vertices are adjacent.
Likewise, the face circles form a packing with circles touching if and only
if the corresponding faces are adjacent. For each edge, there is a pair of
touching vertex circles and a pair of touching face circles. These pairs
touch in the same point, intersecting each other orthogonally.

This circle pattern is unique up to M{\"o}bius transformations.
\end{theorem}

The first published statement and proof of this theorem seems to be contained
in~\cite{BriSch93}. For generalizations, see~\cite{Sch92}, \cite{Riv96},
and~\cite{BS02}, the latter also for a variational proof
(see also Section \ref{sec:patterns} of this article).

Now, mark the centers of the circles with white dots and mark the
intersection points, where two touching pairs of circles intersect each other
orthogonally, with black dots. Draw edges from the center of each circle to
the intersection points on its periphery. You obtain a quad-graph with
bicolored vertices. Since, furthermore, the black vertices have degree four,
the white vertices may be labeled $\sph$ and $\cir$ to make the quad-graph an
S-quad-graph.


Now let us construct the spheres intersecting $S^2$ orthogonally along the
circles marked by $\sph$. Connecting the centers of touching spheres, one
obtains a {\em Koebe polyhedron}: a convex polyhedron with all edges tangent
to the sphere $S^2$. Moreover, the circles marked with $\cir$ are inscribed
into the faces of the polyhedron, see Fig.~\ref{fig:Koebe-polyhedron}. Thus
we have a polyhedral discrete S-isothermic surface. The discrete S-isothermic
surface is given by the spheres $\sph$ and the circles $\cir$.

Thus, Theorem~\ref{thm:orthoKoebe} implies the following theorem.

\begin{theorem}
\label{thm:polyKoebe}
Every polytopal cell decomposition of the
sphere can be realized by a polyhedron with edges tangent to the
sphere. This realization is unique up to projective
transformations which fix the sphere.

There is a simultaneous realization of the dual polyhedron, such
that corresponding edges of the dual and the original polyhedron
touch the sphere in the same points and intersect orthogonally.
\end{theorem}
The last statement of the theorem follows from the construction if
one interchanges the $\cir$ and $\sph$ labels.

\section{Discrete minimal surfaces}
\label{sec:discrete_minimal}

The following theorem about continuous minimal surfaces is due to
Christoffel~\cite{Chr1867}. For a modern treatment, see~\cite{HJ03}. This
theorem is the basis for our definition of discrete minimal surfaces. We
provide a short proof for the reader's convenience.

\begin{theorem}[Christoffel]\label{thm:christoffel}
  Minimal surfaces are isothermic. An isothermic immersion is a minimal
  surface, if and and only if the dual immersion is contained in a sphere. In
  that case the dual immersion is in fact the
  Gauss map of the minimal surface, up to scale and translation.
\end{theorem}

\begin{proof}
  Let $f$ be an isothermic immersion with normal map $N$. Then
  \begin{equation*}
    \langle N_x,f_x\rangle =\lambda^2 k_1\quad\text{and}\quad\langle 
    N_y,f_y\rangle =\lambda^2 k_2,
  \end{equation*}
  where $k_1$ and $k_2$ are the principal curvature functions of $f$ and
  $\lambda=\|f_x\|=\|f_y\|$.  By equations~\eqref{eq:dual}, the dual isothermic
  immersion $f^*$ has normal $N^*=-N$, and
  \begin{equation*}
    \begin{split}
      \langle N^*_x,f^*_x\rangle &=\langle -N_x,\frac{f_x}{\|f_x\|^2}\rangle 
      =- k_1,\\
      \langle N^*_y,f^*_y\rangle &=\langle -N_y,-\frac{f_y}{\|f_y\|^2}\rangle
      =k_2.
    \end{split}
  \end{equation*}
  Its principal curvature functions are therefore
  \begin{equation*}
    k^*_1=-\frac{k_1}{\lambda^2}\quad\text{and}\quad
    k^*_2=\frac{k_2}{\lambda^2}
  \end{equation*}
  Hence $f$ is minimal (this means $k_1=-k_2$) if and only if $f^*$ is
  contained in a sphere ($k^*_1=k^*_2$). In that case, $f^*$ is the Gauss map
  of $f$ (up to scale and translation), because the tangent planes of $f$ and
  $f^*$ at corresponding points are parallel.
\end{proof}

The idea is to define discrete minimal surfaces as S-isothermic surfaces
which are dual to Koebe polyhedra; the latter being a discrete analogue of
conformal parameterizations of the sphere. By
theorem~\ref{thm:minimal_dual_koebe} below, this leads to the following
definition.

\begin{figure}[tb]
  \centering
  \input{discreteminimal.pstex_t}  
  \caption{Condition for discrete minimal surfaces.}
  \label{fig:discrete_minimal_condition}
\end{figure}
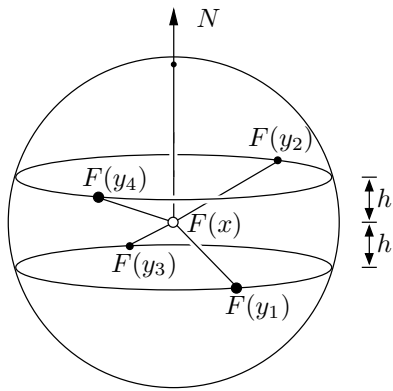

\begin{definition}
  \label{def:discrete_minimal}
  A {\em{}discrete minimal surface}\/ is an S-isothermic discrete surface
  $F:Q\to \R^3$ which satisfies any one of the equivalent conditions
  (i)--(iii) below.
  
  Suppose $x\in Q$ is a white vertex of the quad-graph $Q$ such that $F(x)$
  is the center of a sphere.  Let $y_1\ldots y_{2n}$ be the vertices
  neighboring $x$ in $Q$.  (Generically, $n=2$.) Then $F(y_j)$ are the points
  of contact with the neighboring spheres and simultaneously points of
  intersection with the orthogonal circles. Let $F(y_j)=F(x)+b_j$. (See
  figure~\ref{fig:discrete_minimal_condition}.) Then the following equivalent
  conditions hold:
  \begin{enumerate}
    \renewcommand{\labelenumi}{(\roman{enumi})}
  \item The points $F(x) + (-1)^j b_j$ lie on a circle.
  \item There is an $N\in\R^3$ such that $(-1)^j\,(b_j,N)$ is the same for
    $j=1,\ldots, 2n$. 
  \item There is plane through $F(x)$ such that the points $\{F(y_j)\;|\;j
    \text{ even}\}$ and the points $\{F(y_j)\;|\;j \text{ odd}\}$ lie in
    planes which are parallel to it at the same distance on opposite sides.
  \end{enumerate}
\end{definition}

\begin{remark}\label{rem:maximum_principle}
  The definition implies that a discrete minimal surface is a polyhedral
  surface with the property that every interior vertex lies in the
  convex hull of its neighbors. This is the maximum principle for discrete
  minimal surfaces.
\end{remark}

\begin{examples}
  Fig.~\ref{fig:enneper_catenoid} {\em(left)}\/ shows a discrete minimal
  Enneper surface. Only the circles are shown. A variant of the discrete
  minimal Enneper surface is shown in Fig.~\ref{fig:menneper}. Here, only
  the touching spheres are shown. Fig.~\ref{fig:enneper_catenoid}
  {\em(right)}\/ shows a discrete minimal catenoid. Both spheres and circles
  are shown.  Fig.~\ref{fig:schwarz_scherk} shows a discrete minimal
  Schwarz P-surface and a discrete minimal Scherk tower.
\end{examples}

  These examples are discussed in detail in section~\ref{sec:examples}.

\begin{theorem}\label{thm:minimal_dual_koebe}
  An S-isothermic discrete surface is a discrete minimal surface, if and only
  if the dual S-isothermic surface corresponds to a Koebe polyhedron.
\end{theorem}

\begin{proof}
  That the S-isothermic dual of a Koebe polyhedron is a discrete minimal
  surface is fairly obvious. On the other hand, let $F:Q\to \R^3$ be
  a discrete minimal surface and let $x\in Q$ and $y_1\ldots y_{2n}\in Q$ be
  as in Definition~\ref{def:discrete_minimal}.  Let
  $\widetilde{F}:Q\to \R^3$ be the dual S-isothermic surface. We need
  to show that all circles of $\widetilde{F}$ lie in one and the same sphere
  $S$ and that all the spheres of $\widetilde{F}$ intersect $S$ orthogonally.
  It follows immediately from Definition~\ref{def:discrete_minimal} that the
  points $\widetilde{F}(y_1)\ldots \widetilde{F}(y_{2n})$ lie on a circle
  $c_x$ in a sphere $S_x$ around $\widetilde{F}(x)$. Let $S$ be the sphere
  which intersects $S_x$ orthogonally in $c_x$. The orthogonal circles
  through $\widetilde{F}(y_1)\ldots \widetilde{F}(y_{2n})$ also lie in $S$.
  Hence, all spheres of $\widetilde{F}$ intersect $S$ orthogonally and all
  circles of $\widetilde{F}$ lie in $S$.
\end{proof}

\begin{remark}
  Why do we use S-isothermic surfaces to define discrete minimal surfaces?
  Alternatively, one could define discrete minimal surfaces as the
  surfaces obtained by dualizing discrete (cross-ratio $-1$) isothermic
  surfaces with all quad-graph vertices in a sphere. Indeed, this definition
  was proposed in~\cite{BP96}. However, it turns out that the class of
  discrete isothermic surfaces is too general to lead to a satisfactory
  theory of discrete minimal surfaces. 

  Every way to define the concept of a discrete isothermic immersion imposes
  an accompanied definition of discrete conformal maps. Since a
  conformal map $\R^2\supset D\rightarrow\R^2$ is just an isothermic
  immersion into the plane, discrete conformal maps should be defined as
  discrete isothermic surfaces that lie in a plane.
  Definition~\ref{def:discrete_isothermic} for isothermic surfaces implies
  the following definition for discrete conformal maps: A discrete conformal
  map is a map from a domain of $\Z^2$ to the plane such that all elementary
  quads have cross-ratio $-1$. The so defined discrete conformal maps are too
  flexible. In particular, one can fix one sublattice containing every other
  point and vary the other one, see~\cite{BP99}.

  Definition~\ref{def:S-isothermic} for S-isothermic surfaces, on the other
  hand, leads to discrete conformal maps that are Schramm's ``circle patterns
  with the combinatorics of the square grid''~\cite{Sch97}. This definition
  of discrete conformal maps has many advantages: First, there is Schramm's
  convergence result {\it (ibid)}. Secondly, orthogonal circle patterns have
  the right degree of rigidity. For example, by Theorem~\ref{thm:orthoKoebe},
  two circle patterns that correspond to the same quad-graph decomposition of
  the sphere differ by a M\"obius transformation. One could say: The only
  discrete conformal maps from the sphere to itself are the M\"obius
  transformations. Finally, a conformal map $f: \R^2\supset D\rightarrow\R^2$
  is characterized by the conditions
  \begin{equation}
    \label{eq:conformality conditions}
    |f_x|=|f_y|,\qquad f_x\perp f_y.
  \end{equation}
  To define discrete conformal maps $f:\mathbb{Z}^2\supset D\to \mathbb{C}$,
  it is natural to impose these two conditions on two different sub-lattices
  (white and black) of $\mathbb{Z}^2$, i.e.\ to require that the edges
  meeting at a white vertex have equal length and the edges at a black vertex
  meet orthogonally. Then the elementary quadrilaterals are orthogonal kites,
  and discrete conformal maps are therefore precisely Schramm's orthogonal
  circle patterns.
\end{remark}

\section{A Weierstrass-type representation}
\label{sec:weierstrass}

In the classical theory of minimal surfaces, the Weierstrass representation
allows the construction of an arbitrary minimal surface from holomorphic data
on the underlying Riemann surface.  We will now derive a formula for discrete
minimal surfaces that resembles the Weierstrass representation formula. An
orthogonal circle pattern in the plane plays the role of the holomorphic
data. The discrete Weierstrass representation describes the S-isothermic
minimal surface that is obtained by projecting the pattern stereographically
to the sphere and dualizing the corresponding Koebe polyhedron.

\begin{theorem}[Weierstrass representation]
\label{thm:weierstrass}
Let $Q$ be an S-quad-graph, and let $c:Q\to\C$ be an orthogonal circle
pattern in the plane: For white vertices $x\in Q$, $c(x)$ is the center of
the corresponding circle, and for black vertices $y\in Q$, $c(y)$ is the
corresponding intersection point. The S-isothermic minimal surface
\begin{gather*}
  F:\big\{x\in Q\,\big| \,x\text{ is labelled }\sph\big\} \to\R^3, \\
  F(x) = \text{ the center of the sphere corresponding to $x$}
\end{gather*}
that corresponds to this circle pattern is given by the following formula.
Let $x_1, x_2\in Q$ be two vertices, both labelled $\sph$, that correspond to
touching circles of the pattern, and let $y\in Q$ be the black vertex between
$x_1$ and $x_2$, which corresponds to the point of contact.
The centers $F(x_1)$ and $F(x_2)$ of the corresponding touching spheres
of the S-isothermic minimal surface $F$ satisfy
\begin{multline}
  \label{eq:discrete_weierstrass}
  F(x_2) - F(x_1) = \\
  \pm  \Re\left(\frac{R(x_2) + R(x_1)}{1+\vert p\vert^2}\,
    \frac{\overline{c(x_2)}-\overline{c(x_1)}}{\vert c(x_2)-c(x_1)\vert}  
    \begin{pmatrix}
      1-p^2\\
      i(1+p^2)\\
      2 p
    \end{pmatrix}
  \right),
\end{multline}
where $p=c(y)$ and the radii $R(x_j)$ of the spheres are
\begin{equation}
\label{eq:R}
  R(x_j) = \left\vert \frac{1+\vert c(x_j)\vert^2
      -\vert c(x_j) - p\vert^2}{2 \vert c(x_j) -p\vert}\right\vert
\end{equation}
The sign on the right hand side of equation~\eqref{eq:discrete_weierstrass} depends on
whether the two edges of the quad-graph connecting $x_1$ with $y$ and $y$
with $x_2$ are labelled {\rm `$+$'} or {\rm `$-$'} (see
Figs.~\ref{fig:plusminusedges} and~\ref{fig:OneTwo+S_quad-graph}
{\em(right)}).
\end{theorem}

\begin{proof}
  Let $s:\C\to S^2\subset\R^3$ be the stereographic
  projection
  $$
  s(p) = \frac{1}{1+\vert p\vert^2}
  \begin{pmatrix}
    2 \Re p \\ 2 \Im p \\ \vert p\vert^2-1
  \end{pmatrix}.
  $$
  Its differential is 
  $$
  ds_p(v)= \Re \left( \frac{2\bar v }{(1 +\vert p\vert^2)^2}
    \begin{pmatrix}
      1-p^2 \\ i(1+p^2) \\ 2p 
    \end{pmatrix}
  \right),
  $$
  and
  $$
  \big\|ds_p(v)\big\|=\frac{2\vert v\vert}{1+\vert p\vert^2},
  $$
  where $\|\cdot\|$ denotes the Euclidean norm.
  
  The edge between $F(x_1)$ and $F(x_2)$ of $F$\/ has length $R_1+R_2$
  (this is obvious) and is parallel to $ds_p\big(c(x_2)-c(x_1)\big)$. Indeed,
  this edge is parallel to the corresponding edge of the Koebe polyhedron,
  which, in turn, is tangential to the orthogonal circles in the unit sphere,
  touching in $c(p)$. The pre-images of these circles in the plane touch in
  $p$, and the contact direction is $c(x_2)-c(x_1)$.  Hence,
  equation~\eqref{eq:discrete_weierstrass} follows from
  $$
  F(x_2) - F(x_1) = \pm\big(R(x_2) + R(x_1)\big)
  \frac{ds_p\big(c(x_2)-c(x_1)\big)}{\big\|ds_p\big(c(x_2)-c(x_1)\big)\big\|}
  \,.
  $$
  
\begin{figure}[tbp]%
\centering%
\input{stereo.pstex_t}%
\caption{How to derive equation~\eqref{eq:R}.}%
\label{fig:stereo}%
\end{figure}
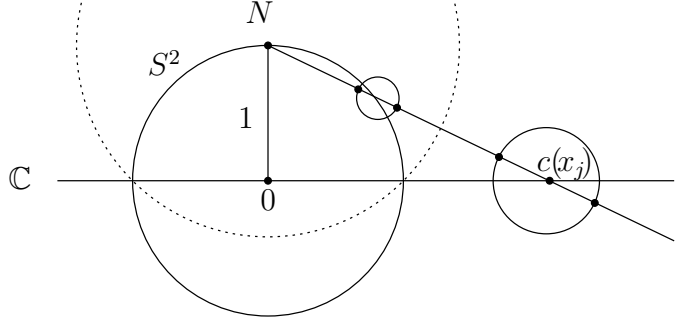

To show equation~\eqref{eq:R}, note that the stereographic projection $s$\/
is the restriction of the reflection on the sphere around the north pole $N$
of $S^2$ with radius $\sqrt{2}$, restricted to the equatorial plane $\C$. See Fig.~\ref{fig:stereo}.  We denote this reflection also by $s$.
Consider the sphere with center $c(x_j)$ and radius $ r=\big|c(x_j)-p\big|, $
which intersects the equatorial plane orthogonally in the circle of the
planar pattern corresponding to $x_j$. This sphere intersects the ray from
the north pole $N$\/ through $c(x_j)$ orthogonally at the distances $d\pm r$
from $N$, where $d$, the distance between $N$ and $c(x_j)$, satisfies $
d^2=1+\big|c(x_j)\big|^2.  $ This sphere is mapped by $s$ to a sphere, which
belongs to the Koebe polyhedron and has radius $1/R_j$.  It intersects the
ray orthogonally at the distances $2/(d\pm r)$. Hence, its radius is
  $$
  1/R_j=\left|\frac{2r}{d^2-r^2}\right|\,.
  $$
  Equation~\eqref{eq:R} follows.
\end{proof}

\section{The associated family}
\label{sec:assoc_family}

Every continuous minimal surfaces comes with an associated family of
isometric minimal surfaces with the same Gauss map. Catenoid and helicoid are
members of the same associated family of minimal surfaces. The concept of an
associated family carries over to discrete minimal surfaces.  In the smooth
case, the members of the associated family remain conformally, but not
isothermically, parameterized. Similarly, in the discrete case, one obtains
discrete surfaces which are not S-isothermic but should be considered as
discrete conformally parameterized minimal surfaces. 

The {\em associated family}\/ of an S-isothermic minimal surface consists of
the one-parameter family of discrete surfaces that are obtained by the
following construction. Before dualizing the Koebe-polyhedron (which would
yield the S-isothermic minimal surface), rotate each edge by an equal angle
in the plane which is tangent to the unit sphere in the point where the edge
touches the unit sphere.

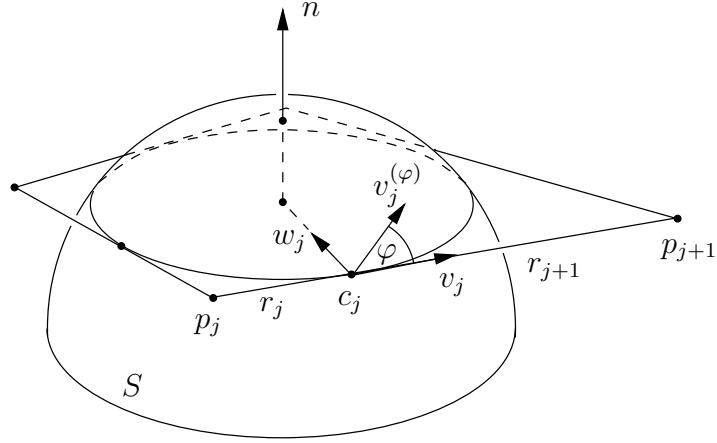
\begin{figure}[tbp]
\begin{center}
\input{assocfig.pstex_t}%
\end{center}
\caption{Proof of Lemma~\ref{lem:assoc_fam}. The vector $v_j^{(\varphi)}$ is
  obtained by rotating $v_j$ in the tangent plane to the sphere at $c_j$.}
\label{fig:assocfig}
\end{figure}

This construction leads to well defined surfaces because of the following
lemma, which is an extension of Lemma~\ref{le:S-dual}. See
Fig.~\ref{fig:assocfig}.

\begin{lemma}
  \label{lem:assoc_fam}
  Let $P$ be a planar polygon with an even number of cyclically ordered edges
  given by the vectors $l_1,\ldots,l_{2n}\in{\R}^3$, $l_1+\ldots+l_{2n}=0$.
  Suppose the polygon has an inscribed circle $c$ with radius $R$, which lies
  in a sphere $S$. Let $r_j$ be the distances from the vertices of $P$ to the
  nearest touching point on the circle:\/ $\| l_j\|=r_j+r_{j+1}$. Rotate each
  vector $l_j$ by an equal angle $\varphi$ in the plane which is tangent to
  $S$ in the point $c_j$ where the edge touches the $S$ to obtain the vectors
  ${l}^{(\varphi)}_1,\ldots,l^{(\varphi)}_{2n}$.  Then the vectors
  $l^{(\varphi)*}_1,\ldots,l^{(\varphi)*}_{2n}$ given by
  $$
  l^{(\varphi)*}_{j}=(-1)^{j}\dfrac{1}{r_j r_{j+1}} l^{(\varphi)}_j
  $$
  satisfy $l^{(\varphi)*}_1+\ldots+l^{(\varphi)*}_{2n}=0$; that is, they
  form a (non-planar) closed polygon.
\end{lemma}

\begin{proof}
  For $j=1,\ldots,2n$ let $(v_j, w_j, n)$ be the orthonormal basis of
  $\R^3$ which is formed by $v_j=l_j/\|l_j\|$, the unit normal $n$ to
  the plane of the polygon $P,$\/ and 
  \begin{equation}
    \label{eq:wj}
    w_j=n\times v_j.
  \end{equation}
  Let $v^{(\varphi)}_j$ be the vector in the tangent plane to the sphere $S$
  at $c_j$ that makes an angle $\varphi$ with $v_j$. Then
  \begin{equation*}
    v^{(\varphi)}_j=\cos\varphi\;v_j + \sin\varphi\cos\theta\;w_j +
    \sin\varphi\sin\theta\;n,
  \end{equation*}
  where $\theta$ is the angle between the tangent plane and the plane of the
  polygon. This angle is the same for all edges. Since 
  \begin{equation*}
    l^{(\varphi)*}_{j}=(-1)^{j}\bigg(\dfrac{1}{r_j}+\frac{1}{r_{j+1}}\bigg)\,
    v^{(\varphi)}_j,
  \end{equation*}
  we have to show that 
  \begin{equation*}
    \sum_{j=1}^{2n} (-1)^{j}\bigg(\dfrac{1}{r_j}+\frac{1}{r_{j+1}}\bigg)\,
    \big(\cos\varphi\;v_j + \sin\varphi\cos\theta\;w_j +
    \sin\varphi\sin\theta\;n\big) = 0.
  \end{equation*}
  By Lemma~\ref{le:S-dual},
  \begin{equation*}
    \sum_{j=1}^{2n} (-1)^{j}\bigg(\dfrac{1}{r_j}+\frac{1}{r_{j+1}}\bigg)
    \,v_j  = 0.
  \end{equation*}
  Due to~\eqref{eq:wj},
  \begin{equation*}
    \sum_{j=1}^{2n} (-1)^{j}\bigg(\dfrac{1}{r_j}+\frac{1}{r_{j+1}}\bigg)
     \,w_j = 0 
  \end{equation*}
  as well. Finally, 
  \begin{equation*}
    \sum_{j=1}^{2n} (-1)^{j}\bigg(\dfrac{1}{r_j}+\frac{1}{r_{j+1}}\bigg)
    \,n = 0,
  \end{equation*}
  because it is a telescopic sum.
\end{proof}

The following two theorems are easy to prove. First, the Weierstrass-type
formula of Theorem~\ref{thm:weierstrass} may be extended to the associate
family.

\begin{theorem}
\label{thm:associatedFamily}
  Using the notation of Theorem~\ref{thm:weierstrass}, the discrete surfaces
  $F_{\varphi}$ of the associated family satisfy
  \begin{multline*}
    F_{\varphi}(x_2) - F_{\varphi}(x_1) = \\
    \pm  \Re\left(
      e^{i\varphi}\,\frac{R(x_2) + R(x_1)}{1+\vert p\vert^2}\,
      \frac{\overline{c(x_2)}-\overline{c(x_1)}}{\vert c(x_2)-c(x_1)\vert}
      \begin{pmatrix}
        1-p^2\\
        i(1+p^2)\\
        2 p
      \end{pmatrix}
    \right),
  \end{multline*}
\end{theorem}

\begin{figure}[tb]%
\centering%
\includegraphics[width=.5\textwidth, clip=true]{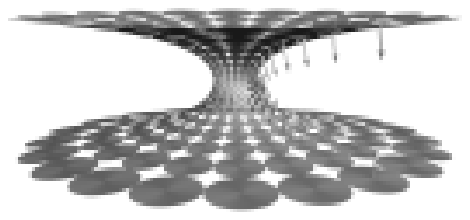}%
\includegraphics[width=.5\textwidth, clip=true]{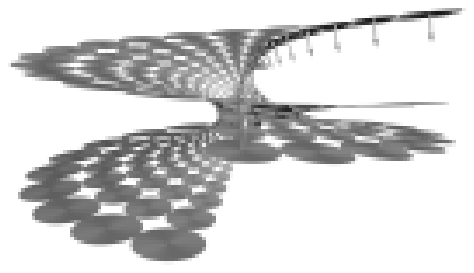}\\%
\includegraphics[width=.5\textwidth, clip=true]{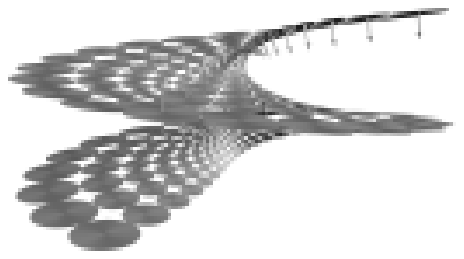}%
\includegraphics[width=.5\textwidth, clip=true]{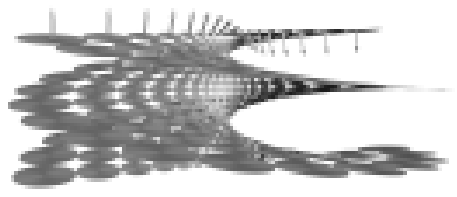}%
\caption{The associated family of the S-isothermic catenoid. The Gauss map is
preserved}%
\label{fig:catenoidAssociatedFamily}%
\end{figure}

Fig.~\ref{fig:catenoidAssociatedFamily} shows the associated family of the
S-isothermic catenoid. 
The essential properties of the associated family of a continuous minimal
surface---that the surfaces are isometric and have the same Gauss
map---carries over to the discrete setting in the following form.

\begin{theorem}
  The surfaces $F_{\varphi}$ of the associated family of an S-isothermic
  minimal surface $F_0$ consist, like $F_0$, of touching spheres.
  The radii of the spheres do not depend on $\varphi$.

  In the generic case, when the quad-graph has $\Z^2$-combinatorics,
  there are also circles through the points of contact, like it is the case
  with $F_0$. The normals of the circles do not depend on $\varphi$.
\end{theorem}

This theorem follows directly from the geometric construction of the
associated family (Lemma~\ref{lem:assoc_fam}).

\section{Convergence}
\label{sec:convergence}

\begin{figure}[tb]
  \centering
  \includegraphics[width=.3333\hsize]{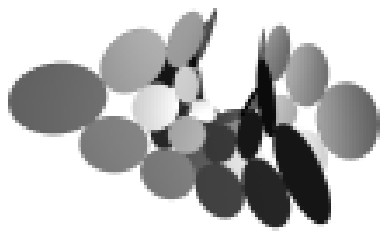}%
  \includegraphics[width=.3333\hsize]{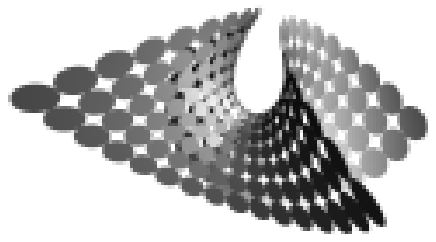}%
  \includegraphics[width=.3333\hsize]{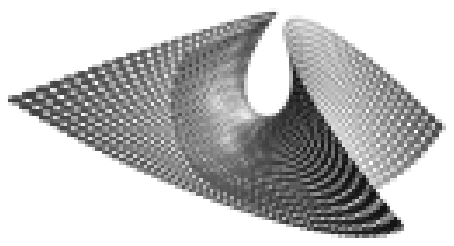}\\

  \caption{A sequence of S-isothermic minimal Enneper surfaces in different discretizations.}
  \label{fig:refinement}
\end{figure}

Schramm has proved the convergence of circle patterns with the combinatorics
of the square grid to meromorphic functions~\cite{Sch97}. Together with the
Weierstrass-type representation formula for S-isothermic minimal surfaces,
this implies the following approximation theorem for discrete minimal
surfaces. Fig.~\ref{fig:refinement} illustrates the convergence of
S-isothermic Enneper surfaces to the continuous Enneper surface.

\begin{theorem}
  Let $D\subset\C$ be a simply connected bounded domain with smooth
  boundary, and let $W\subset\C$ be an open set that contains the
  closure of $D$. Suppose that $F:W\to\R^3$ is a minimal
  immersion without umbilic points in conformal curvature line coordinates.
  There exists a sequence of S-isothermic minimal surfaces
  $\widehat{F}_n:Q_n\to\R^3$ such that the following holds.
  Each $Q_n$ is a simply connected S-quad-graph in $D$ which is a subset of
  the lattice $\frac{1}{n}\Z^2$. If, for $x\in
  D$, we define $\widehat{F}_n(x)$ to be the value of $F_n$ at a point of
  $Q_n$ closest to $x$, then $\widehat{F}_n$ converges to $F$ uniformly with
  error $O(\frac{1}{n})$ on compacts in $D$. In fact, the whole associated
  families $\widehat{F}_{n,\varphi}$ converge to the associated family
  $F_\varphi$ of $F$ uniformly (also in $\varphi$) and with error
  $O(\frac{1}{n})$ on compacts in $D$.
\end{theorem}

\begin{proof}
  Assuming that $F$ is appropriately scaled, 
  \begin{equation}\label{eq:continuous_weierstrass}
    F = \Re\int 
    \begin{pmatrix}
      1-g(z)^2 \\ i\big(1+g(z)^2\big) \\ 2g(z) 
    \end{pmatrix}
    \frac{dz}{g'(z)},
  \end{equation}
  where $g:W\to\C$ is a locally injective meromorphic function. By Schramm's
  results (Theorem 9.1 of~\cite{Sch97} and the remark on p.~387), there
  exists a sequence of orthogonal circle patterns $c_n:Q_n\to\C$
  approximating $g$ and $g'$ uniformly and with error $O(\frac{1}{n})$ on
  compacts in $D$. Define $F_n$ by the Weierstrass
  formula~\eqref{eq:discrete_weierstrass} with data $c=c_n$. Using the
  notation of Theorem~\ref{thm:weierstrass}, one finds
  \begin{equation*}
    \frac{1}{n^2}\big(F_n(x_2)-F_n(x_1)\big) \xrightarrow{n\to\infty}
    \frac{1}{n} 
    \begin{pmatrix}
      1-g(y)^2 \\ i\big(1+g(y)^2\big) \\ 2g(y) 
    \end{pmatrix}
    \frac{1}{g'(y)} + O\Big(\frac{1}{n^2}\Big)
  \end{equation*}
  uniformly on compacts in $D$. Setting $\widehat{F}_n=\frac{1}{n^2}F_n$, the
  convergence claim follows. The same reasoning applies to the whole
  associated family of $F$.
\end{proof}

\section{Orthogonal circle patterns in the sphere}
\label{sec:patterns}

In the simplest cases, like the discrete Enneper surface and the discrete
catenoid (Fig.~\ref{fig:enneper_catenoid}), the construction of the
corresponding circle patterns in the sphere can be achieved by elementary
methods, see Section~\ref{sec:examples}. In general, the problem is not
elementary. Developing methods introduced by
Colin~de~Verdi{\`e}re~\cite{CdV91}, the first and third author have given a
constructive proof of the generalized Koebe theorem, which uses a variational
principle~\cite{BS02}. It also provides a method for the
numerical construction of circle patterns (see also~\cite{Spr03}). An
alternative algorithm was implemented in Stephenson's program
{\tt{}circlepack}~\cite{Ste95}.  It is based on methods developed by
Thurston~\cite{Thu}. The first step in both methods is to transfer the
problem from the sphere to the plane by a stereographic projection. Then the
radii of the circles are calculated. If the radii are known, it is easy to
reconstruct the circle pattern. The radii are determined by a set of
nonlinear equations, and the two methods differ in the way in which these
equations are solved.  Thurston-type methods work by iteratively adjusting
the radius of each circle so that the neighboring circles fit around. The
above mentioned variational method is based on the observation that the
equations for the radii are the equations for a critical point of a convex
function of the radii. The variational method involves minimizing this
function to solve the equations.

\begin{figure}[tb]%
\centering%
\includegraphics[width=0.6\textwidth]{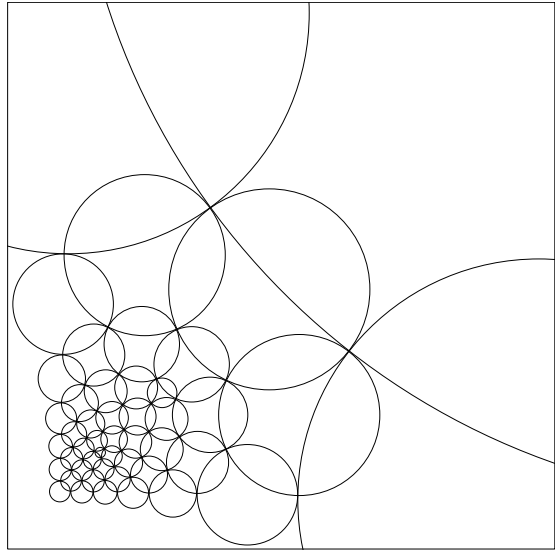}%
\caption{A piece of the circle pattern for a Schwarz P-surface after
  stereographic projection to the plane.}%
\label{fig:euclid_schwarz_pattern}%
\end{figure}
Both of these methods may be used to construct the circle patterns for the
discrete Schwarz P-surface and for the discrete Scherk tower, see
Figs.~\ref{fig:schwarz_scherk} and~\ref{fig:recipe}. One may also take
advantage of the symmetries of the circle patterns and construct only a piece
of it (after stereographic projection) as shown in
Fig.~\ref{fig:euclid_schwarz_pattern}. To this end, one solves the Euclidean
circle pattern problem with Neumann boundary conditions: For boundary
circles, the nominal angle to be covered by the neighboring circles is
prescribed.

However, we have actually constructed the circle patterns for the discrete
Schwarz P-surface and the discrete Scherk tower using an new method suggested
in~\cite{SprPhD}. It is a variational method that works directly on the
sphere. No stereographic projection is necessary; the spherical radii of the
circles are calculated directly. The variational principle for spherical
circle patterns is completely analogous to the variational principles for
Euclidean and hyperbolic patterns presented in \cite{BS02}. We briefly
describe our variational method for circle patterns on the sphere. For a more
detailed exposition, the reader is referred to~\cite{SprPhD}. Here, we will
only treat the case of orthogonally intersecting circles.

The spherical radius $r$ of a non-degenerate circle in the unit sphere
satisfies $0<r<\pi$. Instead of the radii $r$ of the circles, we use the
variables 
\begin{equation}
  \label{eq:rho}
  \rho=\log\tan(r/2). 
\end{equation}
For each circle $j$, we need to find a $\rho_j$ such that the corresponding
radii solve the circle pattern problem. 

\begin{proposition}\label{prop:spherical}
  The radii $r_j$ are the correct radii for the circle pattern if and only if
  the corresponding $\rho_j$ satisfy the following equations, one for each
  circle:

  The equation for circle $j$ is 
  \begin{equation}
    \label{eq:closure}
    2\sum_{\makebox[0pt]{\scriptsize$\text{neighbors } k$}} 
    (\arctan e^{\rho_k-\rho_j} 
    + \arctan e^{\rho_k+\rho_j}) = \Phi_j,
  \end{equation}
  where the sum is taken over all neighboring circles $k$. For each circle $j$,
  $\Phi_j$ is the nominal angle covered by the neighboring circles. It is
  normally $2\pi$ for interior circles, but it differs for circles on the
  boundary or for circles where the pattern branches.  
  
  The equations~\eqref{eq:closure} are the equations for a critical point of
  the functional
  \begin{multline*}
    S(\rho) = \sum_{(j,k)}\big(
    \im\Li_2(ie^{\rho_k-\rho_j}) +
    \im\Li_2(ie^{\rho_j-\rho_k}) \\ -
    \im\Li_2(ie^{\rho_j+\rho_k}) -
    \im\Li_2(ie^{-\rho_j-\rho_k}) - \pi(\rho_j+\rho_k)
    \big) 
    + \sum_j \Phi_j\rho_j.
  \end{multline*}
  Here, the first sum is taken over all pairs $(j,k)$ of neighboring circles,
  the second sum is taken over all circles $j$. The dilogarithm function
  $\Li_2(z)$ is defined by $\Li_2(z)=-\int_0^z \log(1-\zeta)\,d\zeta/\zeta$.
  
  The second derivative of $S(\rho)$ is the quadratic form
  \begin{equation}
    \label{eq:D2S}
    D^2 S =
    \sum_{(j,k)}\Big(\frac{1}{\cosh(\rho_k-\rho_j)}\,(d\rho_k-d\rho_j)^2 -
    \frac{1}{\cosh(\rho_k+\rho_j)}\,(d\rho_k+d\rho_j)^2\Big),
  \end{equation}
  where the sum is taken over pairs of neighboring circles.
\end{proposition}



We provide a proof of Proposition~\ref{prop:spherical} in the Appendix.
Unlike the analogous functionals for euclidean and hyperbolic circle
patterns, the functional $S(\rho)$ is unfortunately not convex : The second
derivative is negative for the tangent vector $v=\sum_j
\partial/\partial\rho_j$. The index is therefore at least $1$. Thus, one
cannot simply minimize $S$ to get to a critical point. However, the following
method seems to work. Define a reduced functional $\widetilde{S}(\rho)$ by
maximizing in the direction $v$:
\begin{equation}
  \label{eq:Stilde}
  \widetilde{S}(\rho)=\max_t S(\rho+t v).
\end{equation}
Obviously, $\widetilde{S}(\rho)$ is invariant under translations in the
direction $v$. Now the idea is to minimize $\widetilde{S}(\rho)$ restricted
to $\sum_j\rho_j=0$. This method has proved to be amazingly
powerful. In particular, it can be used to produce branched circle patterns
in the sphere. It would be very interesting and important to give a
theoretical explanation of this phenomenon.

\section{Constructing discrete minimal surfaces}
\label{sec:constructing}

\begin{figure}[tb]
  \centering
  \parbox[b]{0.4\textwidth}{\includegraphics[width=0.4\textwidth]{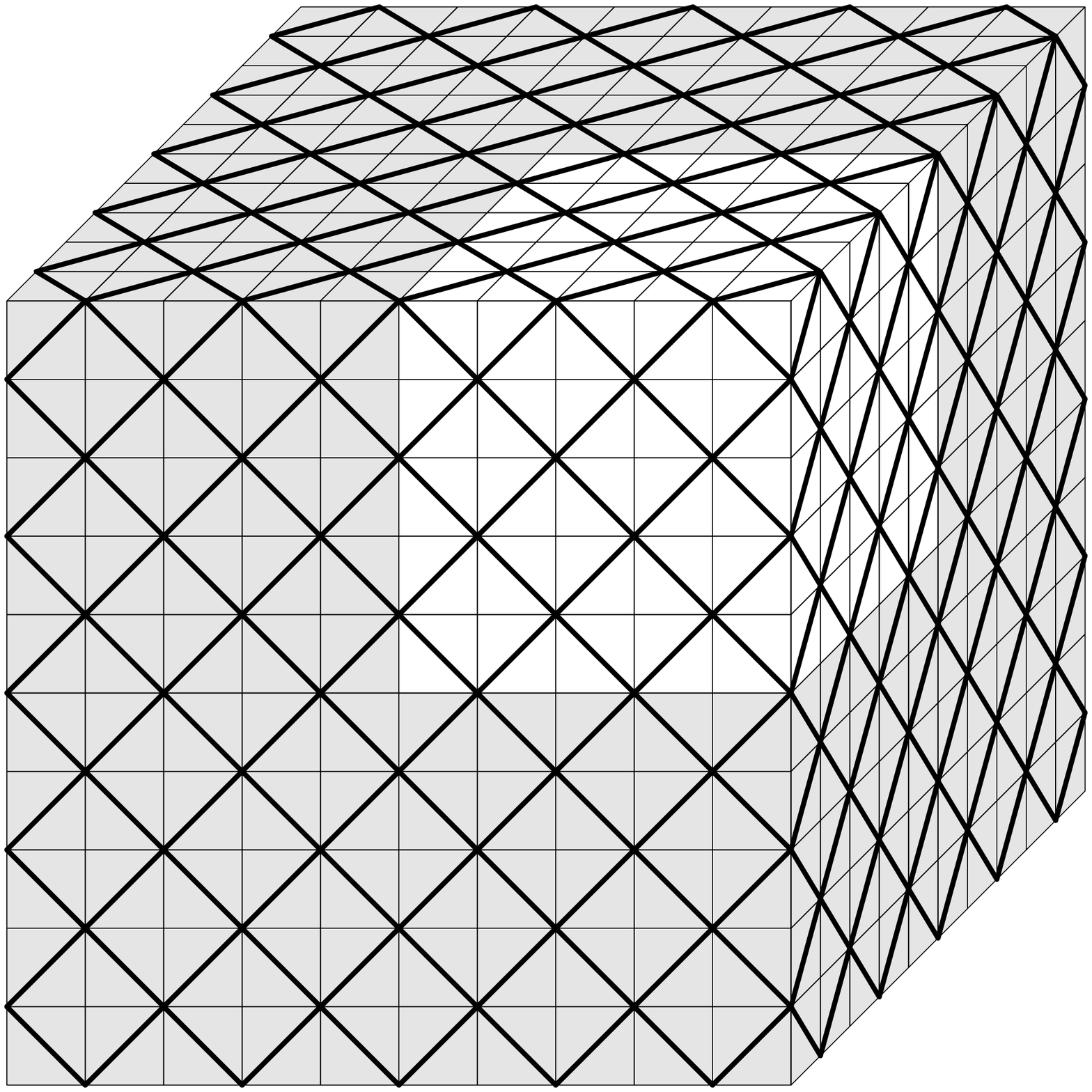}\\
    \centerline{Gauss image of the curvature lines}}
   \raisebox{0.2\textwidth}{\huge\ $\to{}$ }
   \parbox[b]{0.4\textwidth}{\includegraphics[width=0.4\textwidth]{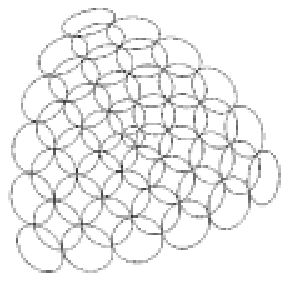}\\
     \centerline{circle pattern}}  \\
   \vspace{\baselineskip}
   \raisebox{0.15\textwidth}{\huge\ $\to{}$ }
   \parbox[b]{0.3\textwidth}{\includegraphics[width=0.3\textwidth]{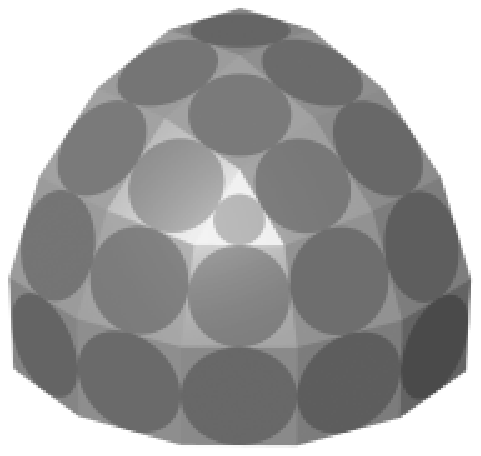}\\
     \centerline{Koebe polyhedron}}
   \raisebox{0.15\textwidth}{\huge\ $\to{}$ }
   \parbox[b]{0.4\textwidth}{\raisebox{0.03\textwidth}{\includegraphics[width=0.4\textwidth]{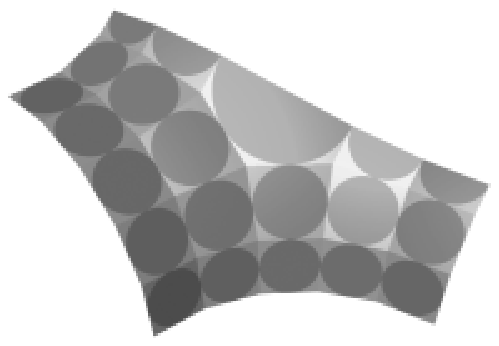}}\\
   \centerline{discrete minimal surface}}
  \caption{Construction of the discrete Schwarz P-surface.}
  \label{fig:recipe}
\end{figure}

Given a specific continuous minimal surface, how does one construct an
analogous discrete minimal surface? In this section we outline the general
method for doing this. By way of an example, Figure~\ref{fig:recipe}
illustrates the construction of the Schwarz P-surface. Details on how to
construct the concrete examples in this paper are explained in
Section~\ref{sec:examples}. The difficult part is finding the right circle
pattern (paragraphs 1 and 2 below). The remaining steps, building the Koebe
polyhedron and dualizing it (Paragraphs 3 and 4) are rather mechanical. 

\paragraph{1. Investigate the Gauss image of the curvature lines.} 
The Gauss map of the continuous minimal surface maps its curvature lines to
the sphere. Obtain a qualitative picture of this image of the curvature lines
under the Gauss map. This yields a quad-graph immersed in the sphere. Here
one has to choose how many curvature lines one wants to use. This corresponds
to a choice of different levels of refinement of the discrete surface. Also,
a choice is made as to which vertices will be black and which will be white.
This choice is usually determined by the nature of the exceptional vertices
corresponding to umbilics and ends. (Umbilics have to be black vertices.)
Only the combinatorics of this quad graph matter.  (Figure~\ref{fig:recipe},
{\em{}top left}.) Generically, the (interior) vertices have degree 4.
Exceptional vertices correspond to ends and umbilic points of the continuous
minimal surface. In the Figure~\ref{fig:recipe}, the corners of the cube are
exceptional.  They correspond to the umbilic points of the Schwarz P-surface.
The exceptional vertices may need to be treated specially. For details see
section~\ref{sec:examples}.

\paragraph{2. Construct the circle pattern.} 
From the quad graph, construct the corresponding circle pattern. White
vertices will correspond to circles, black vertices to intersection points.
Usually, the generalized Koebe theorem is evoked to assert existence and
M{\"o}bius uniqueness of the pattern. The problem of practically calculating
the circle pattern was discussed in Section~\ref{sec:patterns}. Use symmetries
of the surface or special points where you know the direction of the normal
to eliminate the M\"obius ambiguity of the circle pattern.

\paragraph{3. Construct the Koebe polyhedron.} 
From the circle pattern, construct the Koebe polyhedron. Here, a choice is
made as to which circles will become spheres and which will become circles.
The two choices lead to different discrete surfaces close to each other. Both
are discrete analogues of the continuous minimal surface.

\paragraph{4. Discrete minimal surface.} 
Dualize the Koebe polyhedron to obtain a minimal surface.

\paragraph{}If the function $g(z)$ in the Weierstrass
representation~\eqref{eq:continuous_weierstrass} of the continuous minimal
surface is simple enough, it may be that one can construct an orthogonal
circle pattern that is analogous to (or even approximates) this holomorphic
function explicitly by some other means. For example, this is the case for
the Enneper surfaces and the catenoid (see Section~\ref{sec:constructing}).
In this case one does not use Koebe's theorem to construct the circle
pattern.

\section{Examples}
\label{sec:examples}



We now apply the method outlined in the previous section to construct
concrete examples of discrete minimal surfaces. In the case of Enneper's
surface, the orthogonal circle pattern is trivial. The circle patterns for
the higher order Enneper surfaces and for the catenoid are known circle
pattern analogues of the functions $z^a$  and $e^z$. To construct the circle
patterns for the Schwarz P-surface and the Scherk tower, we use Koebe's
theorem.

\subsection{Enneper's surface}
\label{sec:enneper}

The Weierstrass representation of Enneper's surface in conformal curvature
line coordinates is equation~\eqref{eq:continuous_weierstrass} with $g(z)=z$.
The domain is $\C$, and there are no umbilic points. In the domain, the
curvature lines are the parallels to the real and imaginary axes. The Gauss
map embeds the domain into the sphere. 

\begin{floatingfigure}{.42\textwidth}%
\centering%
\vspace{10pt}%
\includegraphics[width=.42\textwidth]{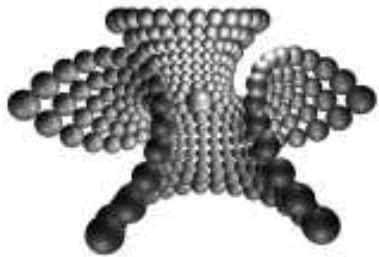}%
\caption{S-isothermic higher order Enneper surface. Only the spheres are 
 shown.}%
\label{fig:menneper}%
\end{floatingfigure}
The quad graph that captures this qualitative behavior of the curvature
lines is the regular square grid decomposition of the plane. There are also
obvious candidates for the circle patterns to use: Take an infinite regular
square grid pattern in the plane. It consists of circles with equal radius
$r$ and centers on a square grid with spacing $\sqrt{2}\,r$. It was shown by
He~\cite{He99} that these patterns are the only embedded and locally finite
orthogonal circle patterns with this quad graph. Project it stereographically
to the sphere and build the Koebe polyhedron. Dualize to obtain a discrete
version of Enneper's surface. See Figs.~\ref{fig:enneper_catenoid}
{\em(left)} and~\ref{fig:refinement}.


\subsection{The higher order Enneper surfaces}
\label{sec:menneper}

As the next example, consider the higher order Enneper
surfaces~\cite{DHKW92}.  Their Weierstrass representation has $g(z)=z^a$.
One may think of them as Enneper surfaces with an umbilic point in the
center.

An orthogonal circle pattern analogue of the maps $z^a$ was introduced
in~\cite{BP99}.  Sectors of these circle patterns were proven to be
embedded~\cite{AB00}, \cite{Aga03}. Stereographic projection to the sphere
followed by dualization leads to S-isothermic analogues of the higher order
Enneper surfaces.  An S-isothermic higher order Enneper surface with a simple
umbilic point ($a=4/3$) is shown in Fig.~\ref{fig:menneper}.

\subsection{The catenoid}
\label{sec:catenoid}

The next most simple example is a discrete version of the catenoid. Here,
$g(z)=e^z$. The corresponding circle pattern in the plane is the
$\operatorname{S-Exp}$ pattern~\cite{BP99}, a discretization of the
exponential map. The underlying quad-graph is $\Z^2$, with circles
corresponding to points $(m,n)$ with $m+n \equiv 0 \mod 2$. The centers
$c(m,n)$ and the radii $r(n,m)$ of the circles are
\begin{equation*}
  c(n,m) = e^{\alpha n + i \rho m},\quad
  r(n,m) = \sin(\rho)|c(n,m)|,
\end{equation*}
where
\begin{equation*}
    \rho = \pi/N,\quad
    \alpha = {\mathop{\rm arctanh}\nolimits}\big({\textstyle\frac{1}{2}}|1-e^{2i\rho}|\big).
\end{equation*}
(It is {\em not}\/ true that $c(m,n)$ is an intersection point if $m+n \equiv
1 \mod 2$.)

The corresponding S-isothermic minimal surface is shown in
Fig.~\ref{fig:enneper_catenoid} {\em(right)}.  The associated family of the
discrete catenoid (see Section~\ref{sec:assoc_family}) is shown in
Fig.~\ref{fig:catenoidAssociatedFamily}.

Other discrete versions of the catenoid have been put forward. A discrete
isothermic catenoid is constructed in \cite{BP96}. This construction can be
generalized in such a way that one obtains the discrete S-isothermic catenoid
described above. This works only because the surface is so particularly
simple.  Then, there is also the discrete catenoid constructed in
\cite{PR02}. It is an area minimizing polyhedral surface. This catenoid is
not related to the S-isothermic catenoid.

\subsection{The Schwarz P-surface}
\label{sec:schwarz-P}

The Schwarz P-surface is a triply periodic minimal surface. It is the
symmetric case in a $2$-parameter family of minimal surfaces with $3$
different hole sizes (only the ratios of the hole sizes matter),
see~\cite{DHKW92}. The domain of the Schwarz P-surface, where the translation
periods are factored out, is a Riemann surface of genus $3$. The Gauss map is
a double cover of the sphere with $8$ branch points. The image of the
curvature line pattern under the Gauss map is shown schematically in
Fig.~\ref{fig:recipe} {\em(top left)}, thin lines. It is a refined cube. More
generally, one may consider three different numbers $m$, $n$, and $k$ of
slices in the three directions. The $8$ corners of the cube correspond to the
branch points of the Gauss map. Hence, not $3$ but $6$ edges are incident
with each corner vertex.  The corner vertices are assigned the label $\cir$.
We assume that the numbers $m$, $n$, and $k$ are even, so that the vertices
of the quad graph may be labelled `$\cir$', `$\sph$', and `$\bullet$'
consistently (see Section~\ref{sec:discrete_S-isothermic}).

To invoke Koebe's theorem (in the form of Theorem~\ref{thm:orthoKoebe}),
forget momentarily that we are dealing with a double cover of the sphere.
Koebe's theorem implies the existence and M{\"o}bius-uniqueness of a circle
pattern as shown in Fig.~\ref{fig:recipe} {\em(top right)}. (Only one eighth
of the complete spherical pattern is shown.)\@ The M{\"o}bius ambiguity is
eliminated by imposing octahedral metric symmetry.

Now lift the circle pattern to the branched cover, construct the Koebe
polyhedron and dualize it to obtain the Schwarz P-surface; see
Fig.~\ref{fig:recipe} {\em(bottom row)}. A fundamental piece of the surface
is shown in Fig.~\ref{fig:schwarz_scherk} {\em(left)}.

We summarize these results in a theorem.

\begin{theorem}
  Given three even positive integers $m$, $n$, $k$, there exists a
  corresponding unique (unsymmetric) S-isothermic Schwarz P-surface.
\end{theorem}

Surfaces with the same ratios $m:n:k$ are different discretizations of the
same continuous Schwarz P-surface.  The cases with $m=n=k$ correspond to the
symmetric Schwarz P-surface.

\subsection{The Scherk tower}
\label{sec:scherk}

Finally, consider Scherk's saddle tower, a simply periodic minimal surface,
which is asymptotic to two intersecting planes. There is a $1$-parameter
family, the parameter corresponding to the angle between the asymptotic
planes, see~\cite{DHKW92}.  An S-isothermic minimal Scherk tower is shown in
Fig.~\ref{fig:schwarz_scherk} {\em(right)}.

\begin{figure}[tbp]
\hfill%
\includegraphics[width=0.125\textwidth]{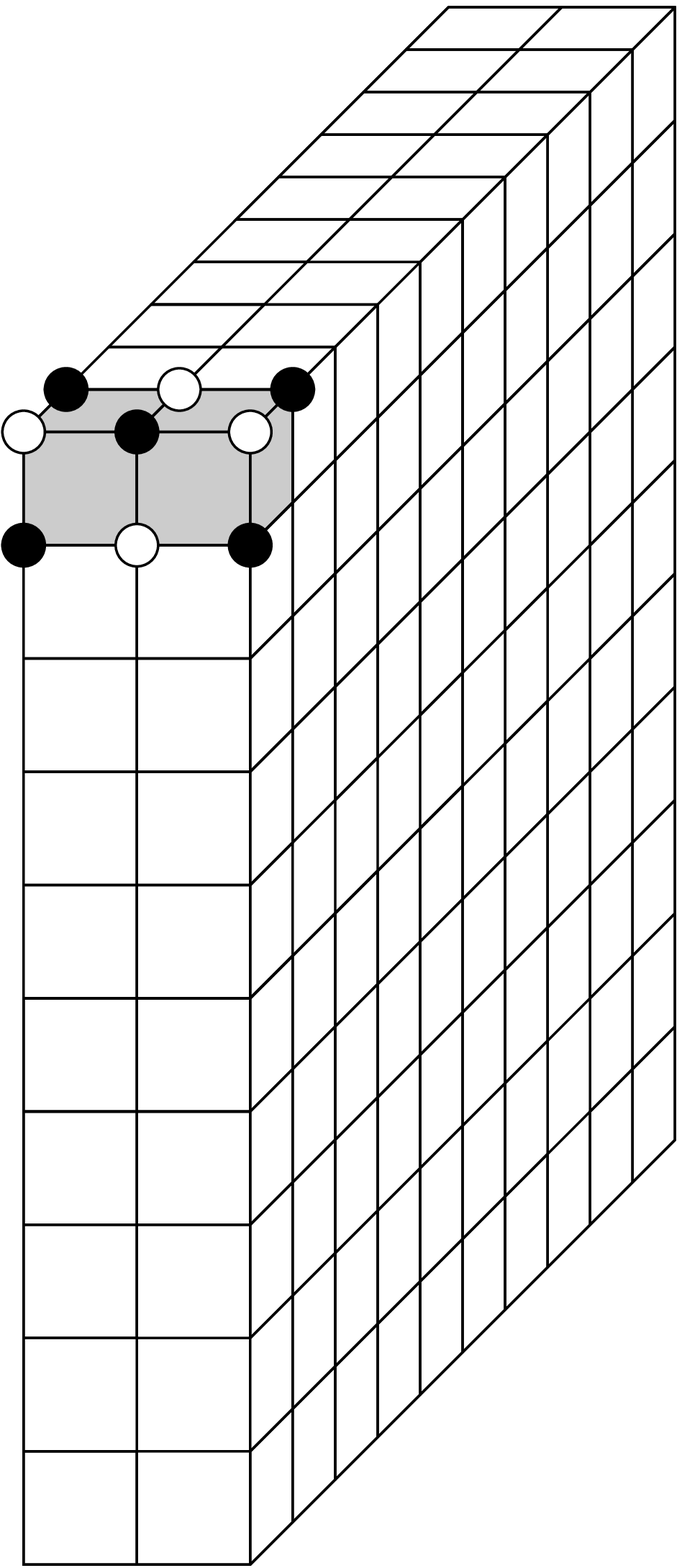}%
\hfill%
\vbox to
0.3\textwidth{\vfill\hbox{\includegraphics[width=0.225\textwidth]{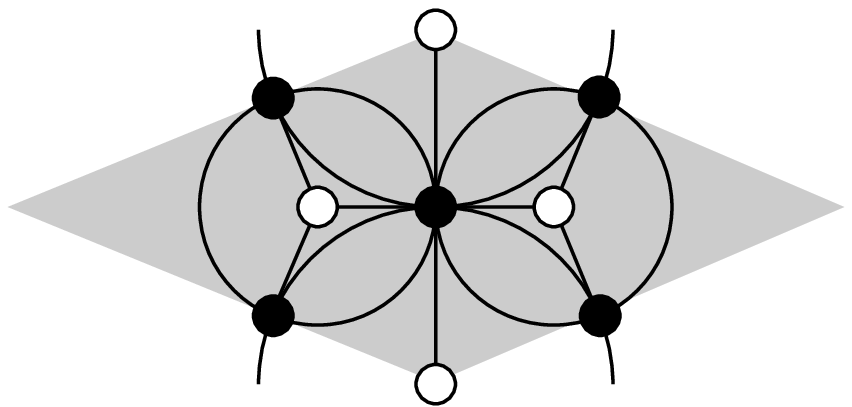}}%
\vfill%
\hbox to 0.225\textwidth{\hfill\bf\Large$\downarrow$\hfill}%
\vfill%
\hbox{\includegraphics[width=0.225\textwidth]{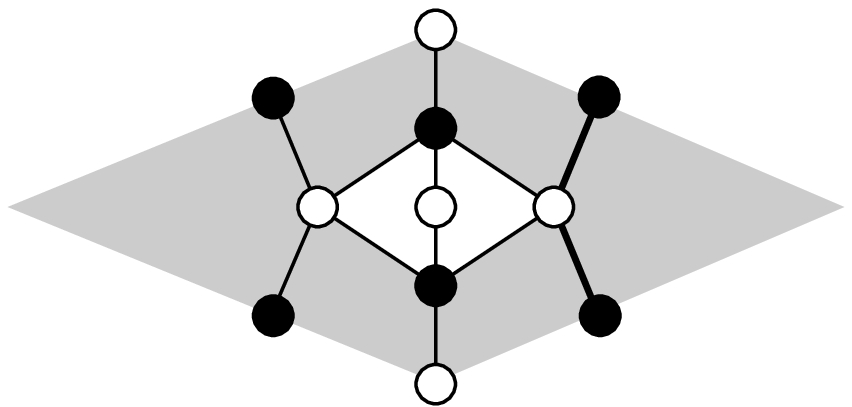}}%
\vfill}%
\hspace*{\fill}
\caption{The combinatorics of the Scherk tower.}
\label{fig:scherkCombinatorics}
\end{figure}

When mapped to the sphere by the Gauss map, the curvature lines of the Scherk
tower form a pattern with four special points, which correspond to the four
half-planar ends. A loop around a special point corresponds to a period of
the surface. In a neighborhood of each special point, the pattern of
curvature lines behaves like the image of the standard coordinate net under
the map $z\mapsto z^2$ around $z=0$. In the discrete setting, the special
points are modeled by pairs of $3$-valent vertices; see
Fig.~\ref{fig:scherkCombinatorics} {\em(left)}. This is motivated by the
discrete version of $z^2$ in~\cite{AB00}. The quad graph we use to construct
the Scherk tower looks like the quad graph for an unsymmetric Schwarz
P-surface with one of the discrete parameters equal to 2. The ratio $m:n$
corresponds to the parameter of the smooth case. Again, by Koebe's theorem,
there exists a corresponding circle pattern, which is made unique by imposing
metric octahedral symmetry.  But now we interpret the special vertices
differently.  Here, they are not branch points. The right hand side of
Fig.~\ref{fig:scherkCombinatorics} shows how they are to be treated: Split
the vertex (and edges) between each pair of $3$-valent vertices in two. Then
introduce new $2$-valent vertices between the doubled vertices. Thus, instead
of pairs of $3$-valent vertices we now have $2$-valent vertices. The newly
inserted edges have length $0$.  Thus, stretching the concept a little bit,
one obtains infinite edges after dualization. This is in line with the fact
that the special points correspond to half-planar ends.

Fig.~\ref{fig:schwarz_scherk} {\em(right)} shows an S-isothermic Scherk
tower.

\begin{theorem}
  Given two even positive integers $m$ and $n$ there exists a
  corresponding unique S-isothermic Scherk tower. 
\end{theorem}

The cases with $m=n$ correspond to the most symmetric Scherk tower, the
asymptotic planes of which intersect orthogonally.

\section*{Appendix. Proof of Proposition~\ref{prop:spherical}}

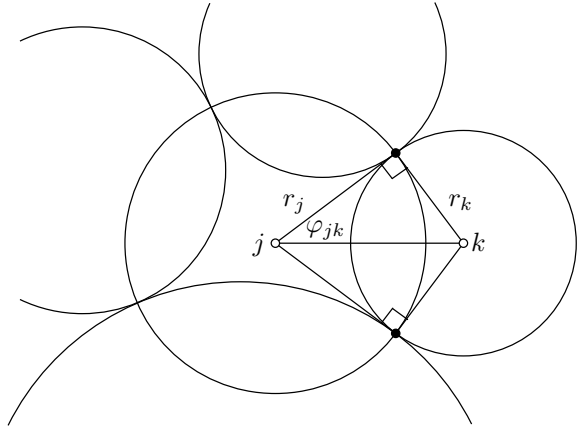
\begin{figure}
\centering
\input{flower.pstex_t}
\caption{A flower of an orthogonal circle pattern.}
\label{fig:flower}
\end{figure}

Figure~\ref{fig:flower} shows a ``flower'' of an orthogonal circle pattern: a
central circle and its orthogonally intersecting neighbors. For simplicity,
it shows a circle pattern in the euclidean plane. We are, however, concerned
with circle patterns in the sphere, where the centers are spherical centers,
the radii are spherical radii and so forth.
The radii of the circles are correct if and only if for each circle the
neighboring circles ``fit around''. This means that for each circle $j$,
\begin{equation*}
  2\sum_{\makebox[0pt]{\scriptsize$\text{neighbors } k$}}\varphi_{jk}
  = \Phi_j,
\end{equation*}
where $\varphi_{jk}$ is half the angle covered by circle~$k$\/ as seen from
the center of circle~$j$, and where normally $\Phi_j=2\pi$\/ except if $j$\/
is a boundary circle or a circle where branching occurs. (In those cases,
$\Phi_j$\/ has some other given value.) Equations~\eqref{eq:closure} follow
from the following spherical trigonometry lemma:

\begin{lemma}
  In a right-angled spherical triangle, let $r_1$ and $r_2$ be the sides
  enclosing the right angle, and let $\varphi$ be the angle opposite side
  $r_2$. Then
  \begin{equation}
    \label{eq:phi}
    \varphi=\arctan e^{\rho_2-\rho_1} + \arctan e^{\rho_2+\rho_1},
  \end{equation}
  where $r$ and $\rho$ are related by equation~\eqref{eq:rho}.
\end{lemma}
\begin{proof}
  Napier's rule says\footnote{In several editions of Bronshtein and
    Semendyayev's {\em Handbook of Mathematics}\/ there is a misprint in the
    corresponding equation.}
  \begin{equation*}
    \tan\varphi=\frac{\tan r_2}{\sin r_1}.
  \end{equation*}
  Equation~\eqref{eq:phi} follows this and the trigonometric
  identity
  \begin{equation}\label{eq:trig_rho}
    \arctan\Big(\frac{\tan r_2}{\sin r_1}\Big)
    =\arctan e^{\rho_2-\rho_1} + \arctan
    e^{\rho_2+\rho_1}.
  \end{equation}
  (To derive equation~\eqref{eq:trig_rho}, start by applying the identity
  \begin{equation*}
    \arctan a + \arctan b = \arctan\frac{a + b}{1 - ab}
  \end{equation*}
  to its right hand side.)
\end{proof}

Now let 
\begin{equation*}
  f(x)=\arctan e^x.
\end{equation*}
Then a primitive function is 
\begin{equation*}
  F(x) = \int_{-\infty}^x f(u)\,du = \im\Li_2(ie^x),
\end{equation*}
(see~\cite{BS02}, \cite{SprPhD}) and the derivative is
\begin{equation*}
  f'(x)=\frac{1}{2\cosh x}\,.
\end{equation*}
Since
\begin{multline}
  \label{eq:S_w_F}
  S(\rho) = \sum_{(j,k)}\big(
  F(\rho_k-\rho_j) + F(\rho_j-\rho_k) - F(\rho_j+\rho_k) - F(-\rho_j-\rho_k)
  - \pi(\rho_j+\rho_k)\big) \\
  + \sum_j \Phi_j\rho_j,
\end{multline}
one obtains after some manipulations that
\begin{equation*}
  \frac{\partial S(\rho)}{\partial\rho_j}=
  - 2\sum_{\makebox[0pt]{\scriptsize$\text{neighbors } k$}} 
    (\arctan e^{\rho_k-\rho_j} 
    + \arctan e^{\rho_k+\rho_j}) + \Phi_j.
\end{equation*}
This proves that Equations~\eqref{eq:closure} are the equations for a
critical point of $S(\rho)$.

Equation~\eqref{eq:D2S} for the second derivative of $S$ is obtained by
taking the second derivative term by term in the first sum of
equation~\eqref{eq:S_w_F}. For example, the second derivative of
$F(\rho_k-\rho_j)$ is $f'(\rho_k-\rho_j)(d\rho_k-d\rho_j)^2$.

This concludes the proof of Proposition~\ref{prop:spherical}.

\vspace{\baselineskip}

\small

\noindent
{\bf Alexander I.~Bobenko}\\
{\bf Tim Hoffmann}\\
{\bf Boris A.~Springborn}

\smallskip\noindent
Technische Universit{\"a}t Berlin \\
Fakult\"at II - Institut f{\"u}r Mathematik\\
Strasse des 17. Juni 136\\
10623 Berlin, Germany 

\smallskip\noindent
\url{bobenko@math.tu-berlin.de}\\
\url{hoffmann@math.tu-berlin.de}\\
\url{springb@math.tu-berlin.de}

\end{document}

%% file: conformalsquare.pstex_t
\begin{picture}(0,0)%
\includegraphics{conformalsquare.pstex}%
\end{picture}%
\setlength{\unitlength}{4144sp}%
\begingroup\makeatletter\ifx\SetFigFont\undefined%
\gdef\SetFigFont#1#2#3#4#5{%
  \reset@font\fontsize{#1}{#2pt}%
  \fontfamily{#3}\fontseries{#4}\fontshape{#5}%
  \selectfont}%
\fi\endgroup%
\begin{picture}(1104,1094)(529,-788)
\put(677,-577){\makebox(0,0)[lb]{\smash{{\SetFigFont{10}{12.0}{\familydefault}{\mddefault}{\updefault}{\color[rgb]{0,0,0}$a$}%
}}}}
\put(665, -7){\makebox(0,0)[lb]{\smash{{\SetFigFont{10}{12.0}{\familydefault}{\mddefault}{\updefault}{\color[rgb]{0,0,0}$b'$}%
}}}}
\put(1327,-626){\makebox(0,0)[lb]{\smash{{\SetFigFont{10}{12.0}{\familydefault}{\mddefault}{\updefault}{\color[rgb]{0,0,0}$b$}%
}}}}
\put(1320, 35){\makebox(0,0)[lb]{\smash{{\SetFigFont{10}{12.0}{\familydefault}{\mddefault}{\updefault}{\color[rgb]{0,0,0}$a'$}%
}}}}
\put(721,-286){\makebox(0,0)[lb]{\smash{{\SetFigFont{10}{12.0}{\familydefault}{\mddefault}{\updefault}{\color[rgb]{0,0,0}$\frac{\mbox{$aa'$}}{\mbox{$bb'$}}=-1$}%
}}}}
\end{picture}%

%% file: sisothermic.pstex_t
\begin{picture}(0,0)%
\includegraphics{sisothermic.pstex}%
\end{picture}%
\setlength{\unitlength}{4144sp}%
\begingroup\makeatletter\ifx\SetFigFont\undefined%
\gdef\SetFigFont#1#2#3#4#5{%
  \reset@font\fontsize{#1}{#2pt}%
  \fontfamily{#3}\fontseries{#4}\fontshape{#5}%
  \selectfont}%
\fi\endgroup%
\begin{picture}(4641,2174)(-97,-2318)
\put(3376,-781){\makebox(0,0)[lb]{\smash{{\SetFigFont{10}{12.0}{\familydefault}{\mddefault}{\updefault}{\color[rgb]{0,0,0}touching spheres}%
}}}}
\put(3376,-1096){\makebox(0,0)[lb]{\smash{{\SetFigFont{10}{12.0}{\familydefault}{\mddefault}{\updefault}{\color[rgb]{0,0,0}orthogonal circles}%
}}}}
\put(3376,-1726){\makebox(0,0)[lb]{\smash{{\SetFigFont{10}{12.0}{\familydefault}{\mddefault}{\updefault}{\color[rgb]{0,0,0}orthogonal kite}%
}}}}
\put(3376,-1411){\makebox(0,0)[lb]{\smash{{\SetFigFont{10}{12.0}{\familydefault}{\mddefault}{\updefault}{\color[rgb]{0,0,0}planar faces}%
}}}}
\end{picture}%

%% file: discreteminimal.pstex_t
\begin{picture}(0,0)%
\includegraphics{discreteminimal.pstex}%
\end{picture}%
\setlength{\unitlength}{4144sp}%
\begingroup\makeatletter\ifx\SetFigFont\undefined%
\gdef\SetFigFont#1#2#3#4#5{%
  \reset@font\fontsize{#1}{#2pt}%
  \fontfamily{#3}\fontseries{#4}\fontshape{#5}%
  \selectfont}%
\fi\endgroup%
\begin{picture}(2228,2314)(803,-2138)
\put(3016,-1321){\makebox(0,0)[lb]{\smash{{\SetFigFont{10}{12.0}{\familydefault}{\mddefault}{\updefault}{\color[rgb]{0,0,0}$h$}%
}}}}
\put(3016,-1051){\makebox(0,0)[lb]{\smash{{\SetFigFont{10}{12.0}{\familydefault}{\mddefault}{\updefault}{\color[rgb]{0,0,0}$h$}%
}}}}
\put(1936, 29){\makebox(0,0)[lb]{\smash{{\SetFigFont{10}{12.0}{\familydefault}{\mddefault}{\updefault}{\color[rgb]{0,0,0}$N$}%
}}}}
\put(2106,-1691){\makebox(0,0)[lb]{\smash{{\SetFigFont{10}{12.0}{\familydefault}{\mddefault}{\updefault}{\color[rgb]{0,0,0}$F(y_1)$}%
}}}}
\put(1251,-926){\makebox(0,0)[lb]{\smash{{\SetFigFont{10}{12.0}{\familydefault}{\mddefault}{\updefault}{\color[rgb]{0,0,0}$F(y_4)$}%
}}}}
\put(2246,-686){\makebox(0,0)[lb]{\smash{{\SetFigFont{10}{12.0}{\familydefault}{\mddefault}{\updefault}{\color[rgb]{0,0,0}$F(y_2)$}%
}}}}
\put(1426,-1436){\makebox(0,0)[lb]{\smash{{\SetFigFont{10}{12.0}{\familydefault}{\mddefault}{\updefault}{\color[rgb]{0,0,0}$F(y_3)$}%
}}}}
\put(1876,-1201){\makebox(0,0)[lb]{\smash{{\SetFigFont{10}{12.0}{\familydefault}{\mddefault}{\updefault}{\color[rgb]{0,0,0}$F(x)$}%
}}}}
\end{picture}%

%% file: stereo.pstex_t
\begin{picture}(0,0)%
\includegraphics{stereo.pstex}%
\end{picture}%
\setlength{\unitlength}{4144sp}%
\begingroup\makeatletter\ifx\SetFigFont\undefined%
\gdef\SetFigFont#1#2#3#4#5{%
  \reset@font\fontsize{#1}{#2pt}%
  \fontfamily{#3}\fontseries{#4}\fontshape{#5}%
  \selectfont}%
\fi\endgroup%
\begin{picture}(4011,1921)(-668,-1328)
\put(-653,-560){\makebox(0,0)[lb]{\smash{{\SetFigFont{12}{14.4}{\familydefault}{\mddefault}{\updefault}{\color[rgb]{0,0,0}$\mathbb C$}%
}}}}
\put(856,-688){\makebox(0,0)[lb]{\smash{{\SetFigFont{12}{14.4}{\familydefault}{\mddefault}{\updefault}{\color[rgb]{0,0,0}$0$}%
}}}}
\put(2509,-452){\makebox(0,0)[lb]{\smash{{\SetFigFont{12}{14.4}{\familydefault}{\mddefault}{\updefault}{\color[rgb]{0,0,0}$c(\!x_j\!)$}%
}}}}
\put(721,-196){\makebox(0,0)[lb]{\smash{{\SetFigFont{12}{14.4}{\familydefault}{\mddefault}{\updefault}{\color[rgb]{0,0,0}$1$}%
}}}}
\put(766,434){\makebox(0,0)[lb]{\smash{{\SetFigFont{12}{14.4}{\familydefault}{\mddefault}{\updefault}{\color[rgb]{0,0,0}$N$}%
}}}}
\put(181,119){\makebox(0,0)[lb]{\smash{{\SetFigFont{12}{14.4}{\familydefault}{\mddefault}{\updefault}{\color[rgb]{0,0,0}$S^2$}%
}}}}
\end{picture}%

%% file: assocfig.pstex_t
\begin{picture}(0,0)%
\includegraphics{assocfig.pstex}%
\end{picture}%
\setlength{\unitlength}{4144sp}%
\begingroup\makeatletter\ifx\SetFigFont\undefined%
\gdef\SetFigFont#1#2#3#4#5{%
  \reset@font\fontsize{#1}{#2pt}%
  \fontfamily{#3}\fontseries{#4}\fontshape{#5}%
  \selectfont}%
\fi\endgroup%
\begin{picture}(4023,2692)(1639,-3516)
\put(2308,-3260){\makebox(0,0)[lb]{\smash{{\SetFigFont{12}{14.4}{\familydefault}{\mddefault}{\updefault}{\color[rgb]{0,0,0}$S$}%
}}}}
\put(3132,-2748){\makebox(0,0)[lb]{\smash{{\SetFigFont{12}{14.4}{\familydefault}{\mddefault}{\updefault}{\color[rgb]{0,0,0}$r_j$}%
}}}}
\put(3210,-2349){\makebox(0,0)[lb]{\smash{{\SetFigFont{12}{14.4}{\familydefault}{\mddefault}{\updefault}{\color[rgb]{0,0,0}$w_j$}%
}}}}
\put(3386,-983){\makebox(0,0)[lb]{\smash{{\SetFigFont{12}{14.4}{\familydefault}{\mddefault}{\updefault}{\color[rgb]{0,0,0}$n$}%
}}}}
\put(3832,-2435){\makebox(0,0)[lb]{\smash{{\SetFigFont{12}{14.4}{\familydefault}{\mddefault}{\updefault}{\color[rgb]{0,0,0}$\varphi$}%
}}}}
\put(3813,-2056){\makebox(0,0)[lb]{\smash{{\SetFigFont{12}{14.4}{\familydefault}{\mddefault}{\updefault}{\color[rgb]{0,0,0}$v^{(\varphi)}_j$}%
}}}}
\put(4206,-2607){\makebox(0,0)[lb]{\smash{{\SetFigFont{12}{14.4}{\familydefault}{\mddefault}{\updefault}{\color[rgb]{0,0,0}$v_j$}%
}}}}
\put(4732,-2518){\makebox(0,0)[lb]{\smash{{\SetFigFont{12}{14.4}{\familydefault}{\mddefault}{\updefault}{\color[rgb]{0,0,0}$r_{j+1}$}%
}}}}
\put(3601,-2716){\makebox(0,0)[lb]{\smash{{\SetFigFont{12}{14.4}{\familydefault}{\mddefault}{\updefault}{\color[rgb]{0,0,0}$c_j$}%
}}}}
\put(2746,-2851){\makebox(0,0)[lb]{\smash{{\SetFigFont{12}{14.4}{\familydefault}{\mddefault}{\updefault}{\color[rgb]{0,0,0}$p_j$}%
}}}}
\put(5536,-2401){\makebox(0,0)[lb]{\smash{{\SetFigFont{12}{14.4}{\familydefault}{\mddefault}{\updefault}{\color[rgb]{0,0,0}$p_{j+1}$}%
}}}}
\end{picture}%

%% file: flower.pstex_t
\begin{picture}(0,0)%
\includegraphics{flower.pstex}%
\end{picture}%
\setlength{\unitlength}{4144sp}%
\begingroup\makeatletter\ifx\SetFigFont\undefined%
\gdef\SetFigFont#1#2#3#4#5{%
  \reset@font\fontsize{#1}{#2pt}%
  \fontfamily{#3}\fontseries{#4}\fontshape{#5}%
  \selectfont}%
\fi\endgroup%
\begin{picture}(3415,2548)(734,-2510)
\put(2524,-1346){\makebox(0,0)[lb]{\smash{{\SetFigFont{10}{12.0}{\familydefault}{\mddefault}{\updefault}{\color[rgb]{0,0,0}$\varphi_{jk}$}%
}}}}
\put(2386,-1186){\makebox(0,0)[lb]{\smash{{\SetFigFont{10}{12.0}{\familydefault}{\mddefault}{\updefault}{\color[rgb]{0,0,0}$r_j$}%
}}}}
\put(3376,-1186){\makebox(0,0)[lb]{\smash{{\SetFigFont{10}{12.0}{\familydefault}{\mddefault}{\updefault}{\color[rgb]{0,0,0}$r_k$}%
}}}}
\put(3511,-1456){\makebox(0,0)[lb]{\smash{{\SetFigFont{10}{12.0}{\familydefault}{\mddefault}{\updefault}{\color[rgb]{0,0,0}$k$}%
}}}}
\put(2206,-1456){\makebox(0,0)[lb]{\smash{{\SetFigFont{10}{12.0}{\familydefault}{\mddefault}{\updefault}{\color[rgb]{0,0,0}$j$}%
}}}}
\end{picture}%